\documentclass[fleqn]{article}
\usepackage{amssymb,latexsym,theorem}

\newcommand{\tG}{\widetilde{G}}
\newcommand{\tGamma}{\widetilde{\Gamma}}

\newtheorem{prop}{Proposition}[section]
\newtheorem{cor}[prop]{Corollary}
\newtheorem{theo}[prop]{Theorem}
\newtheorem{lemma}[prop]{Lemma}
\newtheorem{res}[prop]{Result} 
\newtheorem{note}{Remark}

\newtheorem{defin}{Definition}
\newtheorem{set}[defin]{Setting}

\title{On two non-building but simply connected compact Tits geometries of type $C_3$}

\author{Antonio Pasini}

\date{}

\begin{document} 

\maketitle

\begin{abstract}
A classification of homogeneous compact Tits geometries of irreducible spherical type, with connected panels and admitting a compact flag-transitive automorphism group acting continuously on the geometry, has been obtained by Kramer and Lytchak \cite{KL1} and \cite{KL2}. According to their main result, all such geometries but two are quotients of buildings. The two exceptions are flat geometries of type $C_3$ and arise from polar actions on the Cayley plane over the division algebra of real octonions. The classification obtained by Kramer and Lytchak does not contain the claim that those two exceptional geometries are simply connected, but this holds true, as proved by Schillewaert and Struyve \cite{SS}. The proof by Schillewaert and Struyve is of topological nature and relies on the main result of \cite{KL1} and \cite{KL2}. In this paper we provide a combinatorial proof of that claim, independent of \cite{KL1} and \cite{KL2}.      
\end{abstract}  

\section{Introduction}

We presume that the reader has some knowledge of diagram geometry, in particular Tits geometries, namely geometries belonging to Coxeter diagrams, and buildings. A celebrated Theorem of Tits \cite{Tits2} states that Tits geometries generally come from buildings. Explicitly: a thick Tits geometry of rank $n \geq 2$  is 2-covered by a building if and only if all of its residues of type $C_3$ are covered by buildings; morover, buildings of rank $n \geq 2$ are 2-simply connected. 

Having mentioned coverings and simple connectedness, I recall that, for $1 \leq k \leq n$, a $k$-{\em covering} of geometries of rank $n$ is a type-preserving morphism which induces isomorphims on rank $k$ residues (with the convention that an $n$-covering is just an isomorphism), the domain of a $k$-covering being called a $k$-{\em cover} of the codomain. A geometry is said to be $k$-{\em simply connected} if it does not admit any proper $k$-cover \cite[Chapter 12]{PasDG}. (It goes without saying that a $k$-covering is {\em proper} if it is not an isomorphism.) I warn that $(n-1)$-coverings are usually called {\em coverings}, for short (which forbids to use the word ``covering" as a free shortening for $k$-covering). Accordingly, a geometry of rank $n$ is said to be {\em simply connected} if it is $(n-1)$-simply connected. In particular, coverings of geometries of rank 3 are 2-coverings  and when we say that a geometry of rank 3 is simply connected we just mean it is 2-simply connected. 

Turning back to the above theorem of Tits, that theorem shows the importance of the investigation of $C_3$ geometries. As noticed by Tits \cite{Tits2}, geometries of type $C_3$ that have no relation at all with buildings can be constructed by some kind of free construction, but more examples exist that are not covered by buildings. Classifying them all is perhaps hopeless. 

Nevertheless, with the help of some reasonable additional hypotheses, something can be done. For instance, the following is well known (Aschabacher \cite{Asch}, Yoshiara \cite{Yoshi}):

\begin{theo}\label{prel1}
There exists a unique flag-transitive finite thick $C_3$-geometry which is not a building. It is simply connected and its automorphism group is isomorphic to the alternating group $\mathrm{Alt}(7)$. 
\end{theo} 

The exceptional geometry of Theorem \ref{prel1} is called the $\mathrm{Alt}(7)$-{\em geometry} (also {\em Neumaier geometry} after its discoverer Neumaier \cite{N}). Calling the elements of a $C_3$ geometry {\em points}, {\em lines} and {\em planes} as explained by the following picture 

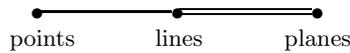
\begin{picture}(310,33)(0,0)
\put(50,18){$\bullet$}
\put(53,21){\line(1,0){50}}
\put(103,18){$\bullet$}
\put(106,20){\line(1,0){50}}
\put(106,22){\line(1,0){50}}
\put(156,18){$\bullet$}
\put(42,8){\small{points}}
\put(97,8){\small{lines}}
\put(146,8){\small{planes}}
\end{picture}

\noindent
the $\mathrm{Alt}(7)$-geometry has $7$ points, $35$ lines and $15$ planes. Moreover, all of its points are incident with all of its planes, namely this geometry is {\em flat}. Referring to \cite{N} (also \cite{Rees} and \cite[\S\S 6.4.2, 12.6.4]{PasDG}) for more details on the $\mathrm{Alt}(7)$ geometry, I only add a few remarks on flat $C_3$-geometries, since we shall deal with them again in this paper. In every flat $C_3$-geometry $\Gamma$ the plane-line system is a linear space $\cal L$, namely a design where any two distinct points (planes of $\Gamma$) belong to a unique common block (line), and the point-line system is a 2-design $\cal D$, possiby with repeated blocks. The geometry $\Gamma$ is obtained by identifying the blocks of $\Delta$ with the lines of $\cal L$ via a suitable bijection.

A number of flag-transitive locally finite (even finite) thick Tits geometries of irreducible type are known that admit the $\mathrm{Alt}(7)$-geometry as a proper residue (see e.g. Buekenhout and Pasini \cite[Section 3]{BP} for a survey), but none of them belongs to a diagram of spherical type. Indeed, as proved by Aschbacher \cite{Asch}, the $\mathrm{Alt}(7)$-geometry cannot occur as a rank 3 residue in any flag-transitive finite thick Tits geometry of irreducible spherical type and rank $n > 3$. Moreover, no finite thick building admits proper quotients (Brouwer and Cohen \cite{BC}). Consequently, 

\begin{cor}\label{prel1bis}
Apart from the $\mathrm{Alt}(7)$-geometry, all flag-transitive finite thick Tits geometries of irreducible spherical type are buildings.
\end{cor} 

Results in the same vein as Theorem \ref{prel1} and Corollary \ref{prel1bis} have recently been obtained by Kramer and Lytchak \cite{KL1} and \cite{KL2} for  compact Tits geometries with connected panels admitting a flag-transitive and compact group of automorhism acting continuously on $\Gamma$.  Before to report on those results, I must explain what a compact geometry is and what we mean when saying that it admits connected panels. 

Let $\Gamma$ be a geometry over a (finite) set of types $I$. Assume that for every $i\in I$ a compact Hausdorff topology is given on the set $\Gamma_i$ of $i$-elements of $\Gamma$ and let ${\cal V}_i$ be the topological space thus defined on $\Gamma_i$. For every $J\subseteq I$ the set $\Gamma_J$ of $J$-flags of $\Gamma$ is a subspace, say ${\cal V}_J$, of the product space $\prod_{j\in J}{\cal V}_j$. If ${\cal V}_J$ is closed (equivalently, compact) for every $J\subseteq I$ then $\Gamma$ is said to be a {\em compact geometry}. (We warn that this definition is not literally the same as in \cite[\S 2.1]{KL1}, but it is equivalent to it.)  When saying that $\Gamma$ has {\em connected panels} we mean that, for every type $i\in I$, the $i$-panels of $\Gamma$ are connected as subspaces of ${\cal V}_i$ (or of ${\cal V}_I$, if we regard panels as sets of chambers). 

With $\Gamma$ a compact geometry as defined above, let $G$ be a flag-transitive group of type-preserving automorphisms of $\Gamma$. Suppose that $G$ is a locally compact topological group (we recall that for topological groups local compactness entails Hausdorff) and that $G$ acts continuously on ${\cal V}_i$ for every $i\in I$ (explicitly, the function from $G\times{\cal V}_i$ to ${\cal V}_i$ that maps $(g,x)\in G\times{\cal V}_i$ onto $g(x)\in{\cal V}_i$ is continuous). Then the pair $(\Gamma, G)$ is called a {\em homogeneous compact geometry} \cite[\S 2.1]{KL1}.
We call $\Gamma$ and $G$ the {\em geometric support} and the {\em group} of $(\Gamma,G)$.  
 
If $(\Gamma, G)$ is a homogeneous compact geometry then $G$ also acts continuously on ${\cal V}_J$ for every $J\subseteq I$. Moreover, for every flag $X\in \Gamma_J$, the stabilizer $G_X$ of $X$ in $G$ is closed in $G$ and the coset space $G/G_X$ is homeomorphic to ${\cal V}_J$, whence Hausdorff and compact. On the other hand, in view of the homeomorphism ${\cal V}_J\approx G/G_X$, the space ${\cal V}_J$ can be recovered from the action and the topology of $G$. So, without assuming any topology on the sets $\Gamma_i$ but still assuming a flag-transitive automorphism group $G$ on $\Gamma$ carrying the structure of a locally compact group with $G_X$ closed and $G/G_X$ compact for every flag $X$ of $\Gamma$, each of the sets $\Gamma_i$ carries a unique compact Hausdorff topology such that $G$ acts continuously on the space ${\cal V}_i$ thus defined on $\Gamma_i$. So, $\Gamma$ is turned into a compact geometry and $(\Gamma,G)$ is a homogeneous compact geometry. 

In this way, as noticed in \cite{KL1}, one can see that all buildings of spherical type associated to semisimple or reductive isotropic algebraic groups defined over local fields are (geometric supports of) homogeneous compact geometries. 

We add one more definition and a few conventions. Given two homogeneous compact geometries $(\widetilde{\Gamma},\widetilde{G})$ and $(\Gamma,G)$ of rank $n \geq 2$ with compact groups $\widetilde{G}$ and $G$, a {\em compact covering} from $(\widetilde{\Gamma},\widetilde{G})$ to $(\Gamma,G)$ is a $2$-covering $\gamma:\widetilde{\Gamma}\rightarrow\Gamma$ such that $\gamma$ is continuous as a mapping from the space $\widetilde{\cal V}$ of elements of $\widetilde{\Gamma}$ to the space $\cal V$ of elements of $\Gamma$, the group $\widetilde{G}$ normalizes the deck group $D$ of $\gamma$ and $\gamma$ induces a continuous isomorphism from the topological group $\widetilde{G}/\widetilde{G}\cap D$ to the topological group $G$. Clearly, $\widetilde{G}\cap D$ is compact. 

The category of homogeneous compact geometries with compact groups and compact coverings as morphisms is named ${\bf HCG}$ in \cite{KL1}. We have defined compact coverings only for homogenous compact geometries with compact groups since these are the objects of $\bf HCG$. According to this restriction, when we say that a given homogeneous compact geometry $(\Gamma,G)$ with $G$ compact is compactly covered by another homogeneous compact geometry $(\widetilde{\Gamma},\widetilde{G})$, it must be understood that $\widetilde{G}$ too is compact.  

We warn the reader that the name ``compact covering"  is not used in \cite{KL1}. We have introduced it with the hope it can remind the reader of the objects of the category $\bf HCG$.

We say that a homogeneous compact geometry is a Tits geometry (in particular, a building) if its geometric support is a Tits geometry (a building). Accordingly, when saying that a homogeneous compact geometry with compact group is compactly covered by a building we mean that it is compactly covered by a homogenous compact geometry the geometric support of which is a building. It goes witout saying that, when speaking of coverings of geometric supports, we mean coverings in the usual `combinatorial' sense, recalled at the beginning of this Introduction. 

More generally, when we say that $(\Gamma,G)$ has some geometric property which neither refers to the topology of $\Gamma$ nor to the group $G$ (such as being a flat $C_3$-geometry, for instance) we mean that the geometric support $\Gamma$ of $(\Gamma, G)$ has that property as a diagram geometry. 

We are now ready to state the main result of Kramer and Lytchak \cite{KL1}, \cite{KL2}.    

\begin{theo}\label{prel2}
Let $(\Gamma,G)$ be a homogeneous compact Tits geometry of type $C_3$ with connected panels and compact group $G$. Then either $(\Gamma,G)$ is compactly covered by a building or it is one of two exceptional flat geometries where $G$ is either $((\mathrm{SU}(3)\times\mathrm{SU}(3))/C_3) >\hspace{-2.4mm}\lhd C_2$ or $\mathrm{SO}(3)\times\mathrm{G}_2$ respectively, in its polar action on the Cayley plane of real octonions. Moreover, the geometric supports of these two exceptional geometries are not covered by any building. 
\end{theo}  

It is convenient to have a name for the two exceptional geometries mentioned in Theorem \ref{prel2}. We shall call them $\mathbb{O}\mathrm{P}^2$-{\em geometries} where $\mathbb{O}$ stands for the octonion algebra over the reals and $\mathbb{O}\mathrm{P}^2$ is the Cayley plane, namely the projective plane over $\mathbb{O}$. 

By exploiting Theorem \ref{prel2}, Kramer and Litchak also obtain the following in \cite{KL1} and \cite{KL2}:    

\begin{cor}\label{prel2bis}
Apart from the two $\mathbb{O}\mathrm{P}^2$-geometries, all homogeneous compact Tits geometries of irreducible spherical type, rank at least $2$, with connected panels and compact group, are compactly covered by buildings.
\end{cor}

The two $\mathbb{O}\mathrm{P}^2$-geometries, or rather the group actions giving rise to them, have been firstly discovered by Podest\`{a} and Thorbergsson \cite{PT} and Gorodski and Kollross \cite{GK}, in the context on an investigation of polar actions of Lie groups on symmetric spaces. A purely algebraic construction of (the geometric supports of) these two geometries is given by Schillewaert and Struyve \cite{SS}. We shall report on that construction in the next section. 

Let $(\Gamma,G)$ be any of the two $\mathbb{O}\mathrm{P}^2$-geometries. The reader should be warned that in the final part of Theorem \ref{prel2} it is not claimed that $\Gamma$ is simply connected. It is only stated that the universal cover $\widetilde{\Gamma}$ of $\Gamma$ is not a building. Thus, in view of rest of the statement of Theorem \ref{prel2}, if $\widetilde{\Gamma} \neq \Gamma$ then either $\widetilde{\Gamma}$ is not the geometric support of any homogeneous compact geometry with compact group or, if it is such, no compact covering exists from that homogeneous compact geometry to $(\Gamma,G)$. So, it is natural to ask if $\Gamma$ is simply connected. The following theorem, due to Schillewaert and Struyve \cite{SS}, answers this question in the affirmative:

\begin{theo}\label{prel3}
The geometric support of either of the two $\mathbb{O}\mathrm{P}^2$-geometries is simply connected.
\end{theo}

The proof that Schillewaert and Struyve give for this theorem is of topological nature. They prove that, if $(\Gamma, G)$ is any of the two $\mathbb{O}\mathrm{P}^2$-geometries, then the universal cover $\widetilde{\Gamma}$ of $\Gamma$ carries a compact Hausdorff topology and $G$ lifts to a compact group $\widetilde{G}\leq \mathrm{Aut}(\widetilde{\Gamma})$, so that $(\widetilde{\Gamma},\widetilde{G})$ is a compact cover of $(\Gamma,G)$. Having proved this, the conclusion follows from Theorem \ref{prel2}: necessarily $\widetilde{\Gamma} = \Gamma$. However Schillewaert and Struyve also collect in \cite{SS} a great deal of information of combinatorial nature on homotopies of closed paths of the two $\mathbb{O}\mathrm{P}^2$-geometries. In this paper we shall exploit that information to arrange a combinatorial proof of Theorem \ref{prel3}, with no use of \cite{KL1} or \cite{KL2}. 

\begin{note}
{\rm As the title of \cite{KL2} makes it clear, an error occurs in \cite{KL1}: the $\mathbb{O}\mathrm{P}^2$-geometry associated to $\mathrm{SO}(3)\times\mathrm{G}_2$ is missing in \cite{KL1}. That gap is filled in \cite{KL2}.}
\end{note}  

\section{The two $\mathbb{O}\mathrm{P}^2$-geometries}

A description of the two $\mathbb{O}P^2$-geometries as coset geometries is given by Kramer and Lytchak in \cite{KL1} (for the geometry with group $G = (\mathrm{SU}(3)\times\mathrm{SU}(3))/C_3 >\hspace{-2.4mm}\lhd C_2$) and in \cite{KL2} (for $G = \mathrm{SO}(3)\times\mathrm{G}_2$). On the other hand, Schillewaert and Struyve \cite{SS} propose a purely algebraic construction for these geometries, which we are going to recall in this section. 

\subsection{Algebraic background} 

Let $\mathbb{A}$ be a division algebra over the field $\mathbb{R}$ of real numbers. It is well known that $\mathbb{A}$ has dimension $1$, $2$, $4$ or $8$ over $\mathbb{R}$. Accordingly, $\mathbb{A}$ is either $\mathbb{R}$ iself or the field $\mathbb{C}$ of complex numbers or the division ring $\mathbb{H}$ or real quaternions or the Cayley-Dyckson algebra $\mathbb{O}$ of real octonions. In any case, $\mathbb{A}$ comes with a {\em norm} $|.|:\mathbb{A}\rightarrow\mathbb{R}$ and a {\em conjugation} $\bar{.}:\mathbb{A}\rightarrow\mathbb{A}$. 

Explicitly, when $\mathbb{A} = \mathbb{R}$ then $|.|$ is the usual absolute value and $\bar{.}$ is the identity; if $\mathbb{A} = \mathbb{C}$ then $|.|$ and $\bar{.}$ are the usual modulus and conjugation. When $\mathbb{A} = \mathbb{H}$ then $\mathbb{A}$ can also be regarded as a right $\mathbb{C}$-vector space with canonical basis $\{1, {\bf j}\}$. The $\mathbb{C}$-span $\mathbb{C} = 1\cdot\mathbb{C}$ of $1$ is a subring of $\mathbb{H}$, ${\bf j}^2 = -1$ and $x{\bf j} = {\bf j}\bar{x}$ for any $x\in\mathbb{C}$. The norm and the conjugation of $\mathbb{H}$ map $x+{\bf j}y$ onto $\sqrt{|x|^2+|y|^2}$ and $\bar{x} - {\bf j}y$ respectively. The conjugation of $\mathbb{H}$ is an involutory anti-automorphism.  Clearly, $\{1, {\bf i}, {\bf j}, {\bf ji}\}$ is a basis of $\mathbb{H}$ over $\mathbb{R}$ (the canonical one), where ${\bf i}$ stands for any of the two square roots of $-1$ in $\mathbb{C}$. 

Finally, $\mathbb{O}$ contains $\mathbb{H}$ as a subring and is generated by $\mathbb{H}$ together with an extra element $\bf k$ such that ${\bf k}^2 = -1$ and 
\begin{equation}\label{eq0}
u{\bf k} ~ = ~ {\bf k}\bar{u}, ~~~ \mbox{for}~ u\in \mathbb{H}
\end{equation}
where $\bar{.}$ denotes the conjugation in $\mathbb{H}$ as defined above. Moreover,
\begin{equation}\label{eq1} 
({\bf k}u)v = {\bf k}(vu) = \bar{v}({\bf k}u), ~~~  ({\bf k}u)({\bf k}v) = -v\bar{u}, ~~~ (\forall u, v \in \mathbb{H}).
\end{equation}
Conditions (\ref{eq1}) imply $(uv){\bf k} = v(u{\bf k}) = v({\bf k}\bar{u})$. Jointly with (\ref{eq0}) they also imply that the elements of $\mathbb{O}$ admit the following representation:  
\begin{equation}\label{eq2}
u + {\bf k}v, ~~~ \mbox{for}~u, v \in \mathbb{H}.
\end{equation}
In spite of (\ref{eq2}), the multiplication of $\mathbb{O}$ does not yield a $\mathbb{H}$-vector space on $\mathbb{O}$, as it follows from the first equality of (\ref{eq1}) and the fact that $\mathbb{H}$ is non-commutative. More precisely, $\mathbb{O}$ does carry a $\mathbb{H}$-vector space structure, as it is clear from (\ref{eq2}), but the scalar multiplication of that space is not the multiplication of $\mathbb{O}$ restricted to $\mathbb{O}\times\mathbb{H}$. On the other hand, for $x, y \in \mathbb{C}$ we have
\[({\bf k}x)y = {\bf k}(yx) = {\bf k}(xy),\]
\[(({\bf k}{\bf j}x)y = ({\bf k}(x{\bf j})y = {\bf k}((y(x{\bf j})) = {\bf k}((yx){\bf j}) = ({\bf k}{\bf j})(yx) = ({\bf k}{\bf j})(xy).\]
So, the multiplication of $\mathbb{O}$ restricted to $\mathbb{O}\times\mathbb{C}$ defines a $4$-dimensional $\mathbb{C}$-vector space on $\mathbb{O}$, with $\{1, {\bf j}, {\bf k}, {\bf k}{\bf j}\}$ as the canonical basis. Needless to say, $\{1, ~{\bf i}, ~{\bf j}, ~ {\bf j}{\bf i}, ~ {\bf k}, ~{\bf k}{\bf i}, ~ {\bf k}{\bf j}, ~{\bf k}({\bf j}{\bf i})\}$ is a basis of $\mathbb{O}$ over $\mathbb{R}$ (the canonical one).  
 
The norm and the conjugation of $\mathbb{O}$ map $u+{\bf k}v$ onto $\sqrt{|u|^2+|v|^2}$ and $\bar{u} - {\bf k}v$ respectively. The conjugation of $\mathbb{O}$ is an involutory anti-automorphism. 

In any case, the norm of $\mathbb{A}$ induces a positive definite $\mathbb{R}$-bilinear form $(.|.)_\mathbb{R}$ which maps $(x,y)\in\mathbb{A}\times\mathbb{A}$ onto the real part $\mathrm{Re}(\bar{x}y)$ of the product $\bar{x}y$. Clearly, $|x| = \sqrt{(x,x)_\mathbb{R}}$. 
We denote by $\perp_\mathbb{R}K$ the orthogonal complement of a subspace $K$ of $\mathbb{A}$ with respect to $(.|.)_\mathbb{R}$.

Let $\mathbb{F}$ be $\mathbb{R}$ or $\mathbb{C}$, with $\mathbb{F} = \mathbb{R}$ when $\mathbb{A} = \mathbb{R}$. Regarded $\mathbb{F}$ as a subfield of $\mathbb{A}$ in the usual way, namely as the $\mathbb{F}$-span of $1$, we set $\mathrm{Pu}_\mathbb{F}(\mathbb{A}) := \perp_\mathbb{R}\mathbb{F}$ (in particular, $\mathrm{Pu}_\mathbb{F}(\mathrm{A}) = 0$ when $\mathbb{A} = \mathbb{F}$). Clearly, $\mathrm{Pu}_\mathbb{F}(\mathbb{A})$ is a subspace of the $\mathbb{F}$-vector space $\mathbb{A}$ and $\mathbb{A} = \mathbb{F}\oplus\mathrm{Pu}_\mathbb{F}(\mathbb{A})$.  The elements of $\mathrm{Pu}_\mathbb{F}(\mathbb{A})$ are said to be $\mathbb{F}$-{\em pure}. 

As $\mathbb{A} = \mathbb{F}\oplus\mathrm{Pu}_\mathbb{F}(\mathbb{A})$, every element $x\in\mathbb{A}$ splits in a unique way as a sum $x = x_1+x_2$ with $x_1\in\mathbb{F}$ and $x_2\in\mathrm{Pu}_\mathbb{F}(\mathbb{A})$. We call $x_1$ and $x_2$ the $\mathbb{F}$-{\em part} and the $\mathbb{F}$-{\em pure part} of $x$. 

When $\mathbb{F} = \mathbb{C}$ we also define a Hermitian inner product $(.|.)_\mathbb{C}:\mathbb{A}\times\mathbb{A}\rightarrow\mathbb{C}$ by 
taking $(x|y)_\mathbb{C}$ equal to the complex part of $\bar{x}y$. Obviously, $\mathrm{Re}((x|y)_\mathbb{C}) = (x|y)_\mathbb{R}$. Hence we also have $|x| = \sqrt{(x|x)_\mathbb{C}}$ for every $x\in \mathbb{A}$. 

The elements of $\mathbb{A}$ of norm $1$ are called {\em unit elements}. Clearly, the set $\mathrm{Un}(\mathbb{A})$ of unit elements of $\mathbb{A}$ is closed under multiplication and taking inverses in $\mathbb{A}$ and
\[\mathbb{A} ~ = ~ \mathrm{Un}(\mathbb{A})\cdot|\mathbb{R}| ~:= ~\{x\cdot|t|~ | ~ x\in \mathrm{Un}(\mathbb{A}), t\in \mathbb{R}\}.\] 
We recall that a homomorphism of $\mathbb{F}$-algebras is an $\mathbb{F}$-linear mapping which also preserves multiplication. In the sequel we shall deal with a particular class of homorphisms of $\mathbb{F}$-algebras, which we shall call sharp $\mathbb{F}$-morphisms. We define them as follows:   

\begin{defin}\label{defin0}
{\rm With $\mathbb{F}$ equal to $\mathbb{R}$ or $\mathbb{C}$, let $\mathbb{A}$ and let $\mathbb{B}$ be two division algebras over $\mathbb{R}$ containing $\mathbb{F}$. When $\mathbb{F} = \mathbb{C}$ both $\mathbb{A}$ and $\mathbb{B}$ can also be regarded as algebras over $\mathbb{C}$. Thus, in any case, both $\mathbb{A}$ and $\mathbb{B}$ are $\mathbb{F}$-algebras. 

A {\em sharp} $\mathbb{F}$-{\em morphism} from $\mathbb{A}$ to $\mathbb{B}$ is a homomorphism of $\mathbb{F}$-algebras from $\mathbb{A}$ to $\mathbb{B}$ which also preserves the inner product $(.|.)_\mathbb{F}$.}   
\end{defin}

Let $\phi:\mathbb{A}\rightarrow\mathbb{B}$ be a sharp $\mathbb{F}$-morphism. Then $\phi$ is injective, since it preserves $(.|.)_\mathbb{F}$. Consequently, $\phi(1) = 1$, hence $\phi(\mathrm{Pu}_\mathbb{F}(\mathbb{A})) \subseteq \mathrm{Pu}_\mathbb{F}(\mathbb{B})$. Moreover $\phi(\mathrm{Un}(\mathbb{A})) \subseteq \mathrm{Un}(\mathbb{B})$. We have $\bar{x} = x^{-1}$ for every unit element $x$. Therefore $\phi(\bar{x}) = \overline{\phi(x)}$ for every $x\in \mathrm{Un}(\mathbb{A})$. Finally, $\phi$ also preserves conjugation.  

As sharp $\mathbb{F}$-morphisms are injective, every sharp $\mathbb{F}$-morphism from $\mathbb{A}$ to $\mathbb{A}$ is an automorphism. We call it a {\em sharp} $\mathbb{F}$-{\em automorphism}. 

\begin{set}\label{setting1}
From now on we assume that $\mathbb{A}$ and $\mathbb{F}$ are as follows: either $\mathbb{A} = \mathbb{H}$ and $\mathbb{F} = \mathbb{R}$ or
$\mathbb{A} = \mathbb{O}$ and $\mathbb{F} = \mathbb{C}$. 
\end{set}

The following is proved in \cite[Proposition 2.1]{SS}:

\begin{lemma}\label{lemma alg1}
With $\mathbb{F}$ and $\mathbb{A}$ as in Setting {\rm \ref{setting1}}, let $a_1, a_2\in \mathrm{Pu}_\mathbb{F}(\mathbb{A})$ and
$b_1, b_2\in \mathrm{Pu}_\mathbb{F}(\mathbb{B})$ be such that $(a_1|a_2)_\mathbb{F} = (b_1|b_2)_\mathbb{F}$, $|a_i| = |b_i|$ for $i = 1, 2$ and $a_1\mathbb{F} \neq a_2\mathbb{F}$. Then there exists a unique sharp $\mathbb{F}$-morphism from $\mathbb{A}$ to $\mathbb{O}$ mapping $a_i$ onto $b_i$ for $i = 1, 2$.  
\end{lemma} 

\begin{lemma}\label{lemma alg2}
Every sharp $\mathbb{R}$-morphism from $\mathbb{H}$ to $\mathbb{O}$ can be estended to a sharp $\mathbb{R}$-automorphism of $\mathbb{O}$.
\end{lemma}
{\bf Proof.} Let $\phi:\mathbb{H}\rightarrow\mathbb{O}$ be a sharp $\mathbb{R}$-morphism. Put ${\bf i}' := \phi({\bf i})$ and ${\bf j}' := \phi({\bf i})$ and recall that $\phi(1) = 1$. Then $\phi(\mathbb{H})$ is the $\mathbb{R}$-span $\mathbb{H}' := \langle 1, {\bf i}', {\bf j}', {\bf j}'{\bf i}'\rangle_\mathbb{R}$ of $\{1, {\bf i}', {\bf j}', {\bf j}'{\bf i}'\}$ and $\phi$ is a sharp $\mathbb{R}$-isomorphism from $\mathbb{H}$ to $\mathbb{H}'$. We can construct a copy $\mathbb{O}'$ of $\mathbb{O}$ starting from $\mathbb{H}'$ instead of $\mathbb{H}$ and, if ${\bf k}'$ is the element of $\mathbb{O}'$ corresponding to $\bf k$, a sharp $\mathbb{R}$-isomorphism $\psi:\mathbb{O}\rightarrow\mathbb{O}'$ is uniquely determined which maps $\bf i$, $\bf j$ and $\bf k$ onto ${\bf i}'$, ${\bf j}'$ and ${\bf k}'$ respectively, whence concides with $\phi$ in $\mathbb{H}$. If we can choose ${\bf k}'\in \mathbb{O}$, then $\psi$ can also be regarded as a sharp $\mathbb{K}$-automorphism of $\mathbb{O}$ and we have done.

So, it remains to prove that we can choose ${\bf k}'\in\mathbb{O}$, namely $\mathbb{O}$ contains an element ${\bf k}'$ orthogonal to $\mathbb{H}$ and such that $({\bf k}')^2 = -1$. But this is obvious. Indeed every unit element orthogonal to $\mathbb{H}$ has this property. The conclusion follows.  \hfill $\Box$ 
   
\subsection{Construction of the geometries}\label{Gamma-F-A} 

With $\mathbb{A}$ and $\mathbb{F}$ as in Setting \ref{setting1}, let $\mathrm{PG}(\mathbb{A})$ be the projetive space of the $\mathbb{F}$-vector space $\mathbb{A}$. For every non-zero vector $x\in \mathbb{A}$, we denote by $[x]$ the corresponding point of $\mathrm{PG}(\mathbb{A})$ and, for every subset $X$ of $\mathbb{A}$ we put $[X] := \{[x]~|~ x\in X\setminus\{0\}\}$. In particular, if $X$ is a subspace of $\mathbb{A}$ then $[X]$ is the corresponding subspace of $\mathrm{PG}(\mathbb{A})$. 

We write $(.|.)$ instead of $(.|.)_\mathbb{F}$ and $\perp$ instead of $\perp_\mathbb{F}$, for short. As usual, $\mathbb{F}^*$ stands for the multiplicative group of $\mathbb{F}$. Following Schillewaert and Struyve \cite{SS}, we construct a $C_3$-geometry $\Gamma_\mathbb{F}(\mathbb{A})$ as follows. 

\begin{defin}\label{defin1}
{\rm The elements (points, lines and planes) of $\Gamma_\mathbb{F}(\mathbb{A})$ are defined as follows:

\begin{itemize}
\item[(A1)] The {\em points} are the points of $[\mathrm{Pu}_\mathbb{F}(\mathbb{A})]$. 
\item[(A2)] The {\em lines} are the sets of pairs $[x,u] := \{(xt,ut)~|~ t\in \mathbb{F}^*\}$ with $x\in \mathrm{Pu}_\mathbb{F}(\mathbb{A})$,  $u\in \mathrm{Pu}_\mathbb{F}(\mathbb{O})$ and $|x| = |u| \neq 0$.
\item[(A3)] The {\em planes} are the sharp $\mathbb{F}$-morphisms $\phi:\mathbb{A}\rightarrow \mathbb{O}$.
\end{itemize}
The {\em incidence relation} is defined as follows:
\begin{itemize}
\item[(B1)] Every point is incident with all planes; 
\item[(B2)] A line $[x,u]$ and a point $[y]$ are declared to be incident when $y\in x^\perp$; 
\item[(B3)] A line $[x,u]$ and a plane $\phi:\mathbb{A}\rightarrow\mathbb{O}$ are incident precisely when $\phi(x) = u$.
\end{itemize}}
\end{defin}

Clearly, the conditions defining point-line and line-plane incidences do not depend on the particular choice of the pair $(x,u) \in [x,u]$.  It is proved in \cite[Proposition 4.3]{SS} that $\Gamma_\mathbb{F}(\mathbb{A})$ is indeed a $C_3$-geometry. According to clause (B1) of Definition \ref{defin1}, this geometry is flat. 

\begin{lemma}\label{lemma geom2}
Both the following hold:

\begin{itemize}
\item[{\rm (1)}] Two lines $[x,u]$ and $[y,v]$ are coplanar if and only if $(x|y) = (u|v)$. If this is the case, then the unique sharp $\mathbb{F}$-morphism $\phi:\mathbb{A}\rightarrow\mathbb{O}$ such that $\phi(x) = u$ and $\phi(y) = v$ (see Lemma {\rm \ref{lemma alg1}}) is the unique plane incident with both $[x,u]$ and $[y,v]$.
\item[{\rm (2)}] If two lines have two distinct points in common then they have the same set of points. 
\end{itemize}
\end{lemma} 
{\bf Proof.} Claim (1) immediately follows from Lemma \ref{lemma alg1} (see also \cite[Lemma 4.2]{SS}). Claim (2) follows from clause (B2) of Definition \ref{defin1} and the fact that $\mathrm{Pu}_\mathbb{F}(\mathbb{A})$ has dimension $3$ over $\mathbb{F}$ (see also \cite[Lemma 5.1]{SS}). \hfill $\Box$  

\bigskip

The set of points of a line $[x,u]$ is the line $x^\perp\cap\mathrm{Pu}_\mathbb{F}(\mathbb{A})$ of $\mathrm{PG}(\mathrm{Pu}_\mathbb{F}(\mathbb{A}))$. We call it the {\em shadow} of $[x,u]$, also a {\em shadow-line}. With this terminology, we can rephrase claim (2) of Lemma \ref{lemma geom2} as follows:

\begin{cor}\label{cor geom3}
The set of points of $\Gamma_\mathbb{F}(\mathbb{A})$ equipped with the shadow lines as lines coincides with the projective plane $\mathrm{PG}(\mathrm{Pu}_\mathbb{F}(\mathbb{A}))$.
\end{cor}  

\subsection{Automorphism groups}\label{Aut groups} 

Let $\mathrm{Aut}_\mathbb{F}(\mathbb{A})$ and $\mathrm{Aut}_\mathbb{F}(\mathbb{O})$ be the groups of sharp $\mathbb{F}$-automorphisms of $\mathbb{A}$ and $\mathbb{O}$. The product $\mathrm{Aut}_\mathbb{F}(\mathbb{A})\times\mathrm{Aut}_\mathbb{F}(\mathbb{O})$ acts on $\Gamma_\mathbb{F}(\mathbb{A})$ as group of automorphisms. Explicitly, given an element 
$(\alpha,\omega)\in\mathrm{Aut}_\mathbb{F}(\mathbb{A})\times\mathrm{Aut}_\mathbb{F}(\mathbb{O})$, 

\begin{itemize}
\item[] $(\alpha,\omega):[x]\rightarrow[\alpha(x)]$ for every point $[x]$ of $\Gamma_\mathbb{F}(\mathbb{A})$;
\item[] $(\alpha,\omega):[x,u]\rightarrow[\alpha(x), \omega(u)]$ for every line $[x,u]$ of $\Gamma_\mathbb{F}(\mathbb{A})$;
\item[] $(\alpha,\omega):\phi\rightarrow\omega\phi\alpha^{-1}$ for every plane $\phi$ of $\Gamma_\mathbb{F}(\mathbb{A})$.
\end{itemize} 
The first questions one may ask is whether this action is faithful and whether all automorphisms of $\Gamma_\mathbb{R}(\mathbb{A})$ arise in these way. Both questions are answered by Schillewaert and Struyve \cite{SS}, but the answers are different according to whether $(\mathbb{F},\mathbb{A}) = (\mathbb{R},\mathbb{H})$ or $(\mathbb{F},\mathbb{A}) = (\mathbb{C},\mathbb{O})$. 

Let $\mathbb{F} = \mathbb{R}$ and $\mathbb{A} = \mathbb{H}$. Then both questions are answered in the affirmative:
\[\mathrm{Aut}(\Gamma_{\mathbb{R}}(\mathbb{H})) ~ = ~  \mathrm{Aut}_\mathbb{R}(\mathbb{H})\times\mathrm{Aut}_\mathbb{R}(\mathbb{O}) ~ = ~ \mathrm{SO}(3)\times \mathrm{G}_2.\]
(Recall that $\mathrm{Aut}_\mathbb{R}(\mathbb{H}) = \mathrm{SO}(3)$ and $\mathrm{Aut}_\mathbb{R}(\mathbb{O}) = \mathrm{G}_2$.) When $\mathbb{F} = \mathbb{C}$ and $\mathbb{A} = \mathbb{O}$ the answer is sligltly different.  Indeed $\mathrm{Aut}_\mathbb{C}(\mathbb{O})\times\mathrm{Aut}_\mathbb{C}(\mathbb{O})$ acts non-faithfully on $\Gamma_{\mathbb{C}}(\mathbb{O})$, with kernel a group $C_3$ of order $3$ contributed by elements $(\zeta, \zeta)$ with $\zeta$ in the center of $\mathrm{SU}(3)$ (recall that $\mathrm{SU}(3) = \mathrm{Aut}_{\mathbb{C}}(\mathbb{O})$). Moreover, the conjugation in $\mathbb{C}$ also induces an automorphism $\gamma$ of $\Gamma_{\mathbb{C}}(\mathbb{O})$ which, being semi-linear as a mapping of $\mathbb{O}\times\mathbb{O}$, does not belong to $\mathrm{Aut}_\mathbb{C}(\mathbb{O})\times \mathrm{Aut}_\mathbb{C}(\mathbb{O})$. All automorphisms of $\Gamma_{\mathbb{C}}(\mathbb{O})$ belong to the group generated by $(\mathrm{Aut}_\mathbb{C}(\mathbb{O})\times \mathrm{Aut}_\mathbb{C}(\mathbb{O}))/C_3$ and $\gamma$. To sum up,
\[\begin{array}{rcll}
\mathrm{Aut}(\Gamma_{\mathbb{C}}(\mathbb{O})) & = &  ((\mathrm{Aut}_\mathbb{C}(\mathbb{O})\times\mathrm{Aut}_\mathbb{C}(\mathbb{O}))/C_3)>\hspace{-2.4mm}\lhd C_2 & =\\
{} & = &  ((\mathrm{SU}(3)\times\mathrm{SU}(3))/C_3)>\hspace{-2.4mm}\lhd C_2. &
\end{array}\]  

\subsection{Recognizing $\Gamma_{\mathbb{F}}(\mathbb{A})$ as an $\mathbb{O}\mathrm{P}^2$-geometry}

Let $\Gamma := \Gamma_\mathbb{F}(\mathbb{A})$ and $G := \mathrm{Aut}(\Gamma)$. As shown by Schillewaert and Struyve \cite[Section 5]{SS}, in either of the two cases that we have considered $(\Gamma, G)$ is a homogeneous compact geometry. They obtain this conclusion by noticing that in either case $G$ is compact and the stabilizers in $G$ of the flags of $\Gamma$ are closed in $G$, but a direct proof is also possible. We shall briefly sketch it here.  

In order to stick to the notation used in the Introduction of this paper, let $\Gamma_1, \Gamma_2$ and $\Gamma_3$ respectively be the set of points, lines and planes of $\Gamma$. In either case each of $\Gamma_1$, $\Gamma_2$ and $\Gamma_3$ can be equipped with a natural compact topology. 

Explicitly, $\Gamma_1 = [\mathrm{Pu}_\mathbb{F}(\mathbb{A})]$ carries the topology of the real projective plane $\mathbb{R}\mathrm{P}^2$ when $(\mathbb{F},\mathbb{A}) = (\mathbb{R},\mathbb{H})$ and the topology of the complex projective plane $\mathbb{C}\mathrm{P}^2$ when $(\mathbb{F},\mathbb{A}) = (\mathbb{R},\mathbb{H})$. Either of these spaces is both Hausdorff and compact. 

When $(\mathbb{F},\mathbb{A}) = (\mathbb{R},\mathbb{H})$ the line-set $\Gamma_2$ carries the topology of the quotient $(\mathbb{ S}^2\times\mathbb{S}^6)/Z$ of the product space $\mathbb{S}^2\times\mathbb{S}^6\subset\mathbb{R}^{10}$ over the center $Z$ of $\mathrm{SL}(\mathbb{R}^{10})$. When $(\mathbb{F},\mathbb{A}) = (\mathbb{C},\mathbb{O})$ then $\Gamma_2$ carries the topology of the quotient $(U\times U)/\Lambda$ where $U := \{x\in\mathbb{C}^3~|~|x| = 1\}$ is the unital  of $\mathbb{C}^3$ and $\Lambda$ is the group of scalar transformations $\lambda\cdot\mathrm{id}$ of $\mathbb{C}^6$ with $|\lambda| = 1$. Again, either of these spaces is Hausdorff and compact.  

When $(\mathbb{F},\mathbb{A}) = (\mathbb{C},\mathbb{O})$ then $\Gamma_3$ carries the same topology as $\mathrm{Aut}_\mathbb{C}(\mathbb{O}) = \mathrm{SU}(3)$, which is (Hausdorff and) compact. Finally, let $(\mathbb{F},\mathbb{A}) = (\mathbb{R},\mathbb{H})$. Then every sharp $\mathbb{R}$-morphism from $\mathbb{H}$ to $\mathbb{O}$ can be regarded as the restriction of a sharp $\mathbb{R}$-automorphism of $\mathbb{O}$ (Lemma \ref{lemma alg2}). Accordingly, the planes of $\Gamma$ naturally correspond to the cosets $\omega H$ of the elementwise stabilizer $H$ of $\mathbb{H}$ in $G := \mathrm{Aut}_\mathbb{R}(\mathbb{O}) = \mathrm{G}_2$. Thus, $\Gamma_3$ can be regaded as a copy of the quotient-space $G/H$, which is still compact and Hausdorff. 

As in the Introduction, let ${\cal V}_1$, ${\cal V}_2$ and ${\cal V}_3$ be the spaces defined on $\Gamma_1$, $\Gamma_2$ and $\Gamma_3$ as above. It is straighforward to check that $\Gamma_{\{i,j\}}$ is closed in ${\cal V}_i\times{\cal V}_j$ for every choice of $1\leq i < j \leq 3$ and the set of chambers $\Gamma_{\{1,2,3\}}$ is closed in ${\cal V}_1\times{\cal V}_2\times{\cal V}_3$. So, $\Gamma$ is a compact geometry. Either of the groups $\mathrm{Aut}(\Gamma_\mathbb{R}(\mathbb{H})) = \mathrm{SO}(3)\times\mathrm{G}_2$ and $\mathrm{Aut}(\Gamma_\mathbb{C}(\mathbb{O})) = ((\mathrm{SU}(3)\times\mathrm{SU}(3))/C_3)>\hspace{-2.4mm}\lhd C_2$ is compact and acts continuously on ${\cal V}_1$, ${\cal V}_2$ and ${\cal V}_3$. 

It remains to show that the group $G$ acts flag-transitively on $\Gamma$. Clearly, in either case $G$ is transitive on the set of point-line fags of $\Gamma$. So, in order to prove flag-transitivity, we only must show that the stabilizer in $G$ of a given point-line flag $([u], [v,x])$ of $\Gamma$ acts transitively on the set of sharp $\mathbb{F}$-morphisms $\phi$ of $\Gamma$ such that $\phi(v) = x$. This follows from Lemma \ref{lemma alg2}. So, 

\begin{res}
The pair $(\Gamma,G)$ is indeed a homogeneous compact geometry.
\end{res} 

As $G$ acts flag-transitively on $\Gamma$, we can recover $\Gamma$ as a coset-geometry from $G$. Comparing flag-stabilizers, it turns out that, when $(\mathbb{F},\mathbb{A}) = (\mathbb{C},\mathbb{O})$, the pair $(\Gamma,G)$ is just the exceptional geometry considered by Kramer and Lytchak in \cite{KL1} (see also Schillewaert and Struyve \cite{SS}). When $(\mathbb{F},\mathbb{A}) = (\mathbb{R},\mathbb{H})$ then $(\Gamma,G)$ is the exceptional geometry of \cite{KL2}. So, 

\begin{res}
The $C_3$-geometries $\Gamma_\mathbb{R}(\mathbb{H})$ and $\Gamma_\mathbb{C}(\mathbb{O})$ are indeed the (geometric supports of the) two $\mathbb{O}\mathrm{P}^2$-geometries. 
\end{res}

\begin{note}
{\rm The two cases of Setting \ref{setting1} correspond to the two cases of \cite{SS} with $\mathbb{B} = \mathbb{O}$. Schillewaert and Struyve \cite{SS} also consider one more case, with $\mathbb{F} = \mathbb{R}$ and $\mathbb{A} = \mathbb{B} = \mathbb{H}$, which leads to a flat $C_3$-geometry which is a quotient of the building associated to the Chevalley group $\mathrm{O}(7,\mathbb{R})$ and admits $\mathrm{SO}(3)\times\mathrm{SO}(3)$ as flag-transitive automorphism group. This geometry also appears in Rees \cite[\S 1.6, (2.2)(ii)]{Rees} as a member of a larger family of flag-transitive flat $C_3$-geometries, obtained as quotients from $\mathrm{O}(7,K)$-buildings, with $K$ any ordered field. Note that the construction used by Rees \cite{Rees} is  primarily geometric. 

This geometry is indeed worth of further investigations, but I have preferred to leave it aside in order to stick to the subject of this paper.}
\end{note}  

\section{A combinatorial proof of Theorem \ref{prel3}} 
 
\subsection{Preliminaries}

We follow \cite{PasDG} for basics on diagram geometry. We recall that, according to \cite{PasDG}, all geometries are residually connected, by definition. In particular, all geometries of rank at least 2 are connected.  

Throughout this subsection $\Gamma$ is a given geometry of rank $n \geq 2$. Recall that $\Gamma$ can be regarded as a simplicial complex, where the vertices are the elements of the geometry and the simplices are the flags. Moreover, with $\{1, 2,..., n\}$ chosen as the type-set of $\Gamma$, the vertices of the complex are marked by positive integers not greater than $n$, according to their type as elements of $\Gamma$. The incidence graph of $\Gamma$ is just the skeleton of the complex $\Gamma$. 

We use the symbol $\sim$ to denote homotopy of paths in the complex $\Gamma$ and $\pi_1(\Gamma)$ for the fundamental group of the complex $\Gamma$. We recall that $\pi_1(\Gamma)$ comes with a distinct vertex chosen as its basis although, as $\Gamma$ is connected, any element of $\Gamma$ can be chosen as such: if we change the basis we accordingly change the group but not its isomorphism type. 

It is well known (see e.g \cite[\S 12.6.1]{PasDG}) that the geometry $\Gamma$ is simply connected (namely $(n-1)$-simply connected) if and only it is simply connected as a complex, namely $\pi_1(\Gamma)$ is trivial. Equivalently, every closed path is null-homotopic.   

\begin{lemma}\label{reduction1}
For $1\leq i < j \leq n$, let $\Gamma_{i,j}$ be the $\{i,j\}$-truncation of $\Gamma$, namely the subgeometry induced by $\Gamma$ on the set of elements of $\Gamma$ of type $i$ or $j$. Then every path of $\Gamma$ starting and ending at $\Gamma_{i,j}$ (in particular, every closed path based at an element of type $i$ or $j$) is homotopic to a path of $\Gamma_{i,j}$.
\end{lemma}
{\bf Proof.} Let $\alpha = (a_0, a_1,..., a_k)$ be a path of $\Gamma$ with $a_0, a_k\in F_{i,j}$. We argue by induction on the length $k$ of $\alpha$. When $k \leq 1$ there is nothing to prove. Let $k = 2$. If $a_1\in \Gamma_{i,j}$ there is nothing to prove as well. Let $a_1\not\in \Gamma_{i,j}$. By the so-called strong connectedness propery \cite[Theorem 1.18]{PasDG}, the intersection $\mathrm{Res}(a_1)\cap \Gamma_{i,j}$ of the residue $\mathrm{Res}(a_1)$ of $a_1$ with $\Gamma_{i,j}$ contains a path 
\[\beta = (b_0 = a_0, b_1,..., b_{h-1}, b_h = a_2)\]
from $a_0$ to $a_2$. We have $(b_{i-1},b_i) \sim (b_{i-1},a_1,b_i)$ for every $i = 1, 2,..., h$, since $\{b_{i-1},a_1,b_i\}$ is a flag. Moreover, $(a_1,b_i,a_1) \sim (a_1)$ for every $i = 1, 2,..., h$. Therefore
\[\beta \sim \gamma := (b_0, a_1, b_1, a_1, b_2, ... , b_{h-1}, a_1, b_h) \sim (b_0, a_1, b_h) = (a_0, a_1, a_2) = \alpha.\]
The claim is proved. Let now $k > 2$. If $a_{k-1}\in \Gamma_{i,j}$ the claim follows by the inductive hypothesis on the subpath $(a_0, a_1,..., a_{k-1})$. Let $a_{k-1}\not\in \Gamma_{i,j}$. If $a_{k-2}\in \Gamma_{i,j}$ then the conclusion follows by the above on the subpath $(a_{k-2},a_{k-1},a_k)$ and the inductive hypothesis on $(a_0, a_1, ... , a_{k-2})$. Let $a_{k-2}\not\in \Gamma_{i,j}$. Pick an element $c\in \mathrm{Res}(a_{k-2},a_{k-1})\cap\Gamma_{i,j}$ and consider the paths
\[\alpha' := (a_0, a_1, ... , a_{k-2}, c), \quad \alpha'' := (c, a_{k-1}, a_k).\]
The path $\alpha'$ has length $k-1$. So, by the inductive hypothesis, a path $\beta'$ exists in $\Gamma_{i,j}$ from $a_0$ to $c$ such that $\beta'\sim \alpha'$. Similarly, as we have already proved the claim for paths of length $2$, a path $\beta''$ exists in $\Gamma_{i,j}$ from $c$ to $a_k$ such that $\beta''\sim\alpha''$. So, $\beta := \beta'\cdot\beta'' \sim \alpha'\cdot\alpha'' = \alpha$ is a path of $\Gamma_{i,j}$ with the required properties.  \hfill $\Box$  

\bigskip

The following lemma is implicit in \cite[Lemma 12.60]{PasDG}. 

\begin{lemma}\label{residues-homotopy}
Given two elements $v$ and $w$ of $\Gamma$, let $\alpha$ and $\beta$ be two paths of $\Gamma$ from $v$ to $w$. If an element $u$ exists in $\Gamma$ such that its residue $\mathrm{Res}(u)$ contains both $\alpha$ and $\beta$, then $\alpha\sim \beta$.
\end{lemma}
{\bf Proof.} Let $\alpha = (a_0, a_1,..., a_k)$ with $a_0 = v$, $a_k = w$ and $\alpha\subseteq \mathrm{Res}(u)$. For every $i = 1, 2,..., k$ put $\alpha_i = (a_{i-1},u,a_i)$. As $(a_{i-1},a_i)\sim (a_{i-1}, u, a_i)$ and $(u, a_i, u) \sim (u)$, we have
\[\alpha \sim \alpha_1\cdot\alpha_2\cdot...\cdot\alpha_k = (a_0, u, a_1, u, a_2,..., a_{k-1}, u, a_k) \sim (a_0, u, a_k).\]
So, $\alpha \sim (a_0, u, a_k) = (v,u,w)$. Similarly, $\beta\sim(v,u,w)$. Therefore $\alpha\sim\beta$.  \hfill $\Box$ 

\subsection{Peculiar properties of $C_3$-geometries}

From now on $\Gamma$ is a geometry of type $C_3$. The integers $1, 2$ and $3$ are taken as types and stand for points, lines and planes respectively.

\begin{defin}\label{primitive path}
{\rm A {\em primitive path} of $\Gamma$ is a closed path $\alpha := (p, L, q, M, p)$ where $p$ and $q$ are points and $L$ and $M$ lines. If $p = q$ or $L = M$ then $\alpha$ is said to be {\em degenerate}.} 
\end{defin} 

Clearly, degenerate primitive paths are null-homotopic. The following is also well known (Tits \cite[Proposition 9]{Tits2}; also \cite[Corollary 7.39]{PasDG}). 

\begin{lemma}\label{C3 buildings}
The geometry $\Gamma$ is a building if and only if all of its primitive paths are degenerate.
\end{lemma}

The proof of the next lemma is implicit in \cite[Section 6.6]{SS}. We make it explicit.

\begin{lemma}\label{reduction2} 
Every closed path of $\Gamma$ based at a point is homotopic to a primitive path.  
\end{lemma}
{\bf Proof}. Let $\alpha$ be a closed path based at a point $p$. In view of Lemma \ref{reduction1}, we may assume that $\alpha$ is contained in $\Gamma_{1,2}$. So, $\alpha = (p_0, L_1, p_1,..., L_k, p_k)$ where $p_0 = p_k = p$ and, for $i = 1,..., k$, $p_i$ is a point and $L_i$ a line. 
We argue by induction on $k$. If $k = 1$ there is  nothing to prove. Let $k > 1$. Suppose firstly that $L_{i-1}$ and $L_i$ are coplanar. Let $\xi$ be the plane on $L_{i-1}$ and $L_i$ and let $M$ be the line of $\mathrm{Res}(\xi)$ through $p_{i-2}$ and $p_i$. Then $(p_{i-2}, L_{i-1},p_{i-1}, L_i, p_i)\sim (p_{i-2}, M, p_i)$ by Lemma \ref{residues-homotopy}. Accordingly, $\alpha \sim \alpha' := (p_0, L_1,..., p_{i-2}, M, p_i,..., L_k, p_k)$. However $\alpha'$, being shorter than $\alpha$, is homotopic to a primitive path, by the inductive hypothesis. Hence $\alpha$ too is homotopic to a primitive path. 

Assume now that $L_{i-1}$ and $L_i$ are never coplanar, for any $i = 2,..., k$. Choose a plane $\xi_2$ on $L_2$. The residue $\mathrm{Res}(p_1)$ of $p_1$ contains a unique line-plane flag $(M_1, \xi_1)$ such that $L_1$ and $M_1$ are incident with $\xi_1$ and $\xi_2$ respectively. Similarly, $\mathrm{Res}(p_2)$ contains a unique line-plane flag $(M_2,\xi_3)$ such that $L_3$ and $M_2$ are incident with $\xi_3$ and $\xi_2$ respectively. 
Let $q$ be the meet-point of $M_1$ and $M_2$ in $\mathrm{Res}(\xi_2)$, let $M_0$ be the line through $p_0$ and $q$ in $\mathrm{Res}(\xi_1)$ and let
$M_3$ be the line through $p_3$ and $q$ in $\mathrm{Res}(\xi_3)$. By Lemma \ref{residues-homotopy} we have the following homotopies:
\[\begin{array}{rcl}
(p_0, L_1, p_1) & \sim & (p_0, M_0, q, M_1, p_1),\\
(p_1, L_2, p_2) & \sim & (p_1, M_1, q, M_2, p_2),\\
(p_2, L_3, p_3) & \sim & (p_2, M_2, q, M_3, p_3).
\end{array}\]
Therefore
\[\begin{array}{l}
(p_0, L_1, p_1, L_2, p_2, L_3, p_3) = (p_0, L_1, p_1)\cdot( p_1, L_2, p_2)\cdot(p_2, L_3, p_3) \sim\\
\sim (p_0, M_0, q, M_1, p_1)\cdot(p_1, M_1, q, M_2, p_2)\cdot(p_2, M_2, q, M_3, p_3) = \\
= (p_0, M_0, q, M_1, p_1, M_1, q, M_2, p_2, M_2, q, M_3, p_3) \sim (p_0, M_0, q, M_3, p_3).
\end{array}\]
Accordingly, $\alpha$ is homotopic to the path, say $\beta$, obtained by replacing the subpath $(p_0, L_1, p_1, L_2, p_2, L_3, p_3)$ of $\alpha$ with 
$(p_0, M_0, q, M_3, p_3)$. The path $\beta$ is shorther than $\alpha$, whence it is homotopic to a primitive path by the inductive hypothesis. As $\alpha\sim \beta$, the same holds for $\alpha$.  \hfill $\Box$ \\

By Lemma \ref{reduction2} we immediately obtain the following: 
 
\begin{cor}\label{reduction3}
The geometry $\Gamma$ is simply connected if and only if all of its primitive paths are null-homotopic.
\end{cor}

Let $\phi_\Gamma:\widetilde{\Gamma}\rightarrow\Gamma$ be the universal covering of $\Gamma$. As $\widetilde{\Gamma}$ is simply connected, all of its closed paths (in particular, all of its primitive paths) are null-homotopic. A closed path of $\Gamma$ is null-homotopic if and only if it lifts through $\phi_\Gamma$ to a closed path of $\widetilde{\Gamma}$. In particular: 

\begin{cor}\label{primitive cor1}
A primitive path of $\Gamma$ is null-homotopic if and only if it is the $\phi_\Gamma$-image of a primitive path of $\widetilde{\Gamma}$. 
\end{cor} 

\begin{cor}\label{primitive cor2}
The geometry $\Gamma$ is covered by a building if and only if none of its non-degenerate primitive paths is null-homotopic.
\end{cor}
{\bf Proof.} Let $\widetilde{\Gamma}$ be a building. Then, by Lemma \ref{C3 buildings}, no non-degenerate primitive path occurs in $\widetilde{\Gamma}$. By Corollary \ref{primitive cor1}, none of the non-degenerate primitive paths of $\Gamma$ can be null-homotopic. On the other hand, let $\widetilde{\Gamma}$ be not a building. Then $\widetilde{\Gamma}$ admits at least one non-degenerate primitive path $\tilde{\alpha}$, necessarily  null-homotopic since $\widetilde{\Gamma}$ is simply connected. Accordingly, $\alpha := \phi_\Gamma(\tilde{\alpha})$ is a null-homotopic non degenerate primitive path of $\Gamma$.  \hfill $\Box$ 

\begin{defin}\label{shift def}
{\rm Let $\alpha = (p, L, q, M, p)$ be a non-degenerate primitive path. Let $N$ be a line on $q$ coplanar with either of $L$ and $M$ and $r$ a point of $N$. The line $N$ is different from either of $L$ and $M$, as $L$ and $M$ are non-copanar. Let $\xi$ be the plane on $N$ and $L$ and let $L'$ be the line of $\xi$ through $p$ and $r$. Similarly, if $\chi$ is the plane on $N$ and $M$, let $M'$ be the line of $\chi$ through $p$ and $r$. Then $(p, L', r, M', p)$ is a primitive path. We denote it by $\sigma_{q\rightarrow r}^N(\alpha)$ and call it the {\em shift} of $\alpha$ {\em from} $q$ {\em to} $r$ {\em along} $N$. We also say that $N$ is {\em admissible} for the path $\alpha$.}
\end{defin}

\begin{lemma}\label{shift}
Let $\alpha = (p, L, q, M, p)$ be a non-degenerate primitive path, $N$ a line admissible for $\alpha$ and $r$ a point of $N$. Then:
\begin{itemize}
\item[$(1)$] We have $\sigma_{q\rightarrow r}^N(\alpha) = \alpha$ if and only if $r = q$.
\item[$(2)$] The shift $\sigma_{q\rightarrow r}^N(\alpha)$ is a non-degenerate primitive path and the line $N$ is admissible for it. 
\item[$(3)$] $\sigma_{r\rightarrow q}^N(\sigma_{q\rightarrow r}^N(\alpha)) = \alpha$.
\item[$(4)$] $\alpha \sim \sigma_{q\rightarrow r}^N(\alpha)$.
\end{itemize}
\end{lemma} 
{\bf Proof.} Claims (1), (2) and (3) are trivial. Claim (4) can be proved as follows:
\[\begin{array}{cccccc}
(p, L, q, M, p) & \sim & (p, \xi, L, q, M, \chi, p) & \sim & (p, \xi, q, \chi, p) & \sim \\
(p, \xi, N, q, N, \chi, p) & \sim & (p, \xi, N, \chi, p) & \sim & (p, \xi, N, r, N, \chi, p) & \sim \\
(p, \xi, r, \chi, p) & \sim & (p, L', \xi, r, \chi, M', p) & \sim & (p, L', r, M', p). & 
\end{array}\] 
(This is essentially the same argoment as used by Schillewaert and Struyve to prove Lemma 6.6 of \cite{SS}.) \hfill $\Box$   

\subsection{Primitive paths in $\mathbb{O}\mathrm{P}^2$-geometries} 

Henceforth $\Gamma = \Gamma_\mathbb{F}(\mathbb{A})$ (see Section \ref{Gamma-F-A}). Recall that $\Gamma$ and the point-line geometry with the same points as $\Gamma$ and the shadow-lines as lines coincides with $\mathrm{PG}(\mathrm{Pu}_\mathbb{F}(\mathbb{A})) \cong \mathrm{PG}(2,\mathbb{F})$ (Corollary \ref{cor geom3}). In particular, two lines of $\Gamma$ either have just one point in common or have exactly the same points. 

\begin{defin}\label{PL-invariant}
{\rm Let $L$ and $M$ be two lines of $\Gamma$ with the same shadow, namely $L = [a, u]$ and $M = [b,v]$ for $a, b\in \mathrm{Pu}_\mathbb{F}(\mathbb{A})$ and $u, v\in \mathrm{Pu}_\mathbb{F}(\mathbb{O})$ with $|a| = |u|\neq 0$, $|b| = |v|\neq 0$ and $[a] = [b]$. Suppose to have chosen the pairs $(a,u)$ and $(b,v)$ in such a way that $a = b$, as we can. Then we put  $(L|M) := (u|v)/|u||v|$. 

Given a primitive path $\alpha = (p,L,q,M,p)$ we put $\ell(\alpha) := (L|M)$ and we call $\ell(\alpha)$ the {\em line-invariant} of $\alpha$. }
\end{defin}

Clearly, $|(L|M)| \leq 1$ by Cauchy-Schwartz inequality, with equality if and only if $u$ and $v$ are proportional. Moreover $(L|M) = 1$ if and only if $L = M$. So, $\ell(\alpha) \neq 1$ whenever $\alpha$ is non-degenerate. 

It is also clear that the hypothesis $a = b$ is necessary for the above definition of $(L|M)$ to make sense. Indeed, without it, only the modulus  $|(u|v)|/|u||v|$ of $(u|v)/|u||v|$ is determined by the pair $L$ and $M$. This also makes it clear that $(L|M)$ can be defined only when $L$ and $M$ have the same shadow. On the other hand, the particular choice of $a$ in the representations $L = [a,u]$ and $M = [a,v]$ is irrelevant. Indeed, if we replace $a$ with $a' = ta$ for some $t\in \mathbb{F}\setminus\{0\}$ then we must also replace $u$ with $u' = tu$ and $v$ with $v' = tv$. Accordingly, $(u'|v')/|u'||v'| = |t|^2(u|v)/|t^2||u||v| = (u|v)/|u||v|$. 

\begin{note}
{\rm Schillewaert and Struyve \cite{SS} call $\ell(\alpha)$ the $PL$-invariant of $\alpha$.}
\end{note} 

\begin{defin}\label{orthogonal path}
{\rm We say that a primitive path $\alpha = (p,L,q,M,p)$ is {\em orthogonal} if $p\perp q$. Assuming that $\alpha$ is non-degenerate but not that it is orthogonal, an {\em orthogonal shift} of $\alpha$ is a shift $\sigma^N_{q\rightarrow r}(\alpha)$ with $p\perp r$.} 
\end{defin}

\begin{lemma}\label{ortho-shift}
Every non-degenerate primitive path $\alpha = (p,L,q,M,p)$ admits orthogonal shifts along every line $N$ admissible for it and, once $N$ has been chosen, the orhogonal shift of $\alpha$ along $N$ is uniquely determined. Moreover, if $\alpha$ is orthogonal, then $\alpha$ is its own orthogonal shift. 
\end{lemma}
{\bf Proof.} As $N$ is coplanar with either of $L$ and $M$, it has at most one point in common with $L$ or $M$. However $N$ contains $q$. Hence it cannot contain $p$. By Corollary \ref{cor geom3}, the line $p^\perp$ of $\mathrm{Pu}_\mathbb{F}(\mathbb{A})$ meets the shadow of $N$ in just one point. (This argument is the same as in the proof of Lemma 6.6 of Schillewaert and Struyve \cite{SS}.) The first part of the lemma is proved. The last claim of the lemma is obvious. \hfill $\Box$ 

\bigskip

Henceforth we denote by $\sigma^N_\perp(\alpha)$ the orthogonal shift of $\alpha$ along a line $N$ admissible for $\alpha$.   

\begin{lemma}\label{ortho-invariant}
Given a non-orthogonal non-degenerate primitive path $\alpha$ of $\Gamma$ and a line $N$ admissible for $\alpha$, let $\beta = \sigma^N_\perp(\alpha)$ be the orthogonal shift of $\alpha$ along $N$ and let $\ell = \ell(\beta)$ be the line-invariant of $\beta$. 

We can always choose the line $N$ in such a way that $\ell \neq -1$.  
\end{lemma}
{\bf Proof.} We must distinguish two cases and two subcases for each of them.\\

\noindent
{\bf Case 1.} $\Gamma = \Gamma_\mathbb{R}(\mathbb{H})$. Modulo automorphism of $\Gamma$, we can always assume that 
\[\begin{array}{lll}
L = [{\bf j},{\bf j}], & & M = [{\bf j}, {\bf i}m_1+{\bf j}m_2], ~~ m_1^2+m_2^2= 1,\\
p = [{\bf i}], & & q = [{\bf i}q_1 + {\bf ji}q_3], ~~ q_1^2 + q_3^2 = 1.
\end{array}\]
So, $\ell(\alpha) = m_2$. Note that $q_1\neq 0$ (otherwise $p\perp q$, while $\alpha$ is non-orthogonal by assumption) and $q_3 \neq 0$ (otherwise $p = q$). Let $N = [b,x]$ be admissible for $\alpha$, where 
\[\begin{array}{l}
b = {\bf i}b_1 + {\bf j}b_2 + {\bf ji}b_3, ~~ b_1^2 + b_2^2 + b_3^2 = 1,\\ 
x = {\bf i}x_1 + {\bf j}x_2 + {\bf ji}x_3 + {\bf k}x_4 + {\bf ki}x_5 + {\bf kj}x_6 + {\bf k}({\bf ji})x_7, ~~ |x| = 1.
\end{array}\]
Modulo automorphisms of $\mathbb{O}$ that leave $\mathbb{H}$ elementwise fixed, we can always assume that 
\[x = {\bf i}x_1 + {\bf j}x_2 + {\bf ji}x_3 + {\bf k}x_4, ~~ (x_1^2 + x_2^2 + x_3^2 + x_4^2 = 1).\]
For $N$ to be admissible for $\alpha$ the following must hold: $({\bf i}q_1+{\bf ji}q_3|b) = 0$ (namely $q$ belongs to $N$) and $({\bf j}|b) = ({\bf j}|x) = ({\bf i}m_1+{\bf j}m_2|x)$ (Lemma \ref{lemma geom2}, claim (1)). Explicitly: 
\begin{equation}\label{caso1 1}
b_1q_1 + b_3q_3 = 0, 
\end{equation}
and $b_2 = x_2 = m_1x_1+m_2x_2$, namely
\begin{equation}\label{caso1 2}
b_2 = x_2, \hspace{5 mm} m_1x_1 = (1-m_2)b_2.
\end{equation}
Let $r = [{\bf i}r_1 + {\bf j}r_2 + {\bf ji}r_3]$ be the unique point of $\{[b],p\}^\perp$. So, $r_1 = 0$, namely $r = [{\bf j}r_2 + {\bf ji}r_3]$, and 
\begin{equation}\label{caso1 3}
b_2r_2 + b_3r_3 = 0.
\end{equation}
Moreover we assume $r_2^2 + r_3^2 = 1$, as we can. We have already noticed that $q_1 \neq 0$. We also have $r_2 \neq 0$, otherwise equations (\ref{caso1 1}) and (\ref{caso1 3}) force $b_1 = b_3 = 0$, hence $b = \pm{\bf j}$, contrary to the fact that $N$ is coplanar with $L$ and $M$. Thus, by (\ref{caso1 1}) and (\ref{caso1 3}) we obtain 
\begin{equation}\label{caso1 4}
b_1 = -b_3q_3q_1^{-1}, \hspace{5 mm} b_2 = - b_3r_3r_2^{-1}.
\end{equation}
These equations show that $b_3 \neq 0$ (otherwise $b = 0$, which is ridiculous). Recalling that $b_1^2 + b_2^2 + b_3^2 = 1$ now we get 
\begin{equation}\label{caso1 5}
b_3  =  \pm\frac{q_1r_2}{\sqrt{q_1^2 + r_2^2 - q_1^2r_2^2}}  =  \pm\frac{q_1r_2}{\sqrt{q_1^2r_3^2 + 1 - r_3^2}} = \pm\frac{q_1r_2}{\sqrt{1-q_3^2r_3^2}}.
\end{equation}
Equation (\ref{caso1 5}) is equivalent to the following
\[r_2 ~ = ~ \pm\frac{b_3}{\sqrt{b_2^2+b_3^2}},\]
which better shows that the point $r$ depends on the choice of the line $N$ but, in view of the sequel,  (\ref{caso1 5}) is more convenient. We shall now consider two subcases: either $m_2 = -1$ or $-1 < m_2 < 1$ (note that $m_2 = 1$ is impossible, since $m_2 = (L|M)$ and $(L|M) \neq 1$ because $L\neq M$). \\

\noindent
{\bf Subcase 1.1.} $m_2 = -1$. Equivalently, $m_1 = 0$. Then $b_2 = x_2 = 0$ by (\ref{caso1 2}), $r_3 = 0$ by (\ref{caso1 4}) and since $b_3 \neq 0$, whence $r_2 = \pm 1$ (as $r_2^2 + r_3^2 = 1$) and  $b_3 = \pm q_1$ by (\ref{caso1 5}). Consequently, $b_1 = \pm q_3$, since $b_1^2 + b_3^2 = q_1^2 + q_3^2 = 1$. Summarizing:
\[\begin{array}{cccccccc}
m_1 & m_2 & r_2 & r_3 & b_1 & b_2 & b_3 & x_2\\
0 & -1 & \pm 1 & 0 & \pm q_3 & 0 &  \pm q_1 & 0.
\end{array}\]
Let now $\xi$ be the plane on $L$ and $N$ and $\chi$ the plane on $M$ and $N$. Then $\xi$ and $\chi$, regarded as  sharp $\mathbb{R}$-morphisms from $\mathbb{H}$ to $\mathbb{O}$, are uniquely determined by the following conditions (Lemma \ref{lemma alg1}): $\xi({\bf j}) = {\bf j}$, $\chi({\bf j}) = {\bf i}m_1 + {\bf j}m_2$ and $\xi(b) = \chi(b) = x$. By entering the above values for $m_1, m_2$ and $x_2$ we get  
\begin{equation}\label{caso11 1}
\xi({\bf j}) = {\bf j}, ~~~ \chi({\bf j}) = -{\bf j}, ~~~ \xi(b) = \chi(b) = {\bf i}x_1 + {\bf ji}x_3 + {\bf k}x_4.
\end{equation} 
Clearly ${\bf i} = {\bf i}(b_1 - {\bf j}b_3)(b_1-{\bf j}b_3)^{-1} = b(b_1+{\bf j}b_3)$. Therefore, and taking (\ref{caso11 1}) into account,
\begin{equation}\label{caso11 2}
\begin{array}{rcl}
\xi({\bf i}) & = & ({\bf i}x_1 + {\bf ji}x_3 + {\bf k}x_4)(b_1 + {\bf j}b_3),\\
\chi({\bf i}) & = & ({\bf i}x_1 + {\bf ji}x_3 + {\bf k}x_4)(b_1 - {\bf j}b_3).
\end{array}
\end{equation} 
Let now $L'$ and $M'$ be the lines through $p$ and $r$ in $\xi$ and $\chi$ respectively. Then $L' = [a,\xi(a)]$ and $M' = [a,\chi(a)]$ where 
$a = {\bf i}a_1 + {\bf j}a_2 + {\bf k}a_3$ is orthogonal with both $p$ and $r$ and we assume $a_1^2+a_2^2+a_3^2 = 1$, as we can. Orthogonality with $p$ and $r$ forces $a_1 = 0 = a_2$. Therefore $a = \pm{\bf ji}$. Accordingly, and recalling (\ref{caso11 2}), 
\begin{equation}\label{caso11 3}
\begin{array}{rcl}
\xi(a) & = & \pm{\bf j}({\bf i}x_1 + {\bf ji}x_3 + {\bf k}x_4)(b_1 + {\bf j}b_3),\\
\chi(a) & = & \overline{+}{\bf j}({\bf i}x_1 + {\bf ji}x_3 + {\bf k}x_4)(b_1 - {\bf j}b_3). 
\end{array}
\end{equation} 
With $\beta = \sigma^N_\perp(\alpha) = (p,L',r,M',p)$ we have $\ell(\beta) = (\xi(a)|\chi(a))$. Equations (\ref{caso11 3}) allow to explicitly compute the inner product $(\xi(a)|\chi(a))$. We obtain:
\begin{equation}\label{caso11 4}
\begin{array}{rcl}
(\xi(a)|\chi(a)) & = &  x_1^2(b_3^2-b_1)^2+x_3^2(b_3^2-b_1^2)+x_4^2(b_3^2-b_1^2) =\\
  & = & (x_1^2+x_3^2+x_4^2)(b_3^2-b_1)^2 = b_3^2-b_1^2 = q_1^2-q_3^2.
\end{array} 
\end{equation}
So, $(\xi(a)|\chi(a)) = q_1^2-q_3^2$. As $q_1, q_3 \neq 0$, we have $-1 < (\xi(a)|\chi(a)) < 1$. \\

\noindent
{\bf Subcase 1.2.} $m_1 \neq 0$, namely $m_2\neq -1$. In this case the second equation of (\ref{caso1 2}) yields
\begin{equation}\label{caso12 1}
x_1 ~ = ~ \frac{1-m_2}{m_1}b_2.
\end{equation} 
The planes $\xi$ and $\chi$ on $L$ and $N$ and on $M$ and $N$ are determined by the following conditions:
\begin{equation}\label{caso12 2}
\begin{array}{l}
\xi({\bf j}) = {\bf j}, \hspace{5 mm} \chi({\bf j}) = {\bf i}m_1 + {\bf j}m_2, \\
\xi(b) = \chi(b)  =  {\bf i}x_1 + {\bf j}x_2 + {\bf ji}x_3 + {\bf k}x_4 =  ({\bf i}\frac{1-m_2}{m_1} + {\bf j})b_2 + {\bf ji}x_3 + {\bf k}x_4.
\end{array}
\end{equation} 
Moreover, $x_3^2 + x_4^2 = 1- ((1-m_2)^2m_1^{-2}+1)b_2^2 = 1-2(1+m_2)^{-1}b_2^2$. Therefore
\begin{equation}\label{caso12 3}
x_3^2 + x_4^2 ~ = ~ 1 - \frac{2}{1+m_2}b_2^2.
\end{equation}
Now ${\bf i} = (b-{\bf j}b_2)(b_1-{\bf j}b_3)^{-1} = (b-{\bf j}b_2)(b_1+{\bf j}b_3)(b_1^2+b_3^2)^{-1}$. Recalling equations (\ref{caso1 4}), we obtain
\begin{equation}\label{caso12 4}
{\bf i} ~ = ~ (b+{\bf j}\frac{r_3}{r_2}b_3)({\bf j}q_1 -q_3)q_1b_3^{-1}.
\end{equation} 
As in Subcase 1.1, let $L' = [a,\xi(a)]$ and $M' = [a,\chi(a)]$ be the lines through $p$ and $r$ in $\xi$ and $\chi$ respectively, where 
$a = {\bf i}a_1 + {\bf j}a_2 + {\bf k}a_3$ with $|a| = 1$. The vector $a$ is orthogonal with both $p$ and $r$. Orthogonality with $p$ still forces  $a_1 = 0$ but orthogonality with $r$ only implies $a_2r_2 + a_3r_3 = 0$. So $a_2 = - a_3r_3r_2^{-1}$ and the condition $|a| = 1$ implies $a_3 = \pm r_2$. Hence $a_2 = \pm r_3$. Summarizing
\begin{equation}\label{caso12 5}
a = \pm({\bf j}r_3 + {\bf ji}r_2).
\end{equation}
Exploiting (\ref{caso12 2}), (\ref{caso12 4}) and (\ref{caso12 5}), we can compute $\xi(a)$ and $\chi(a)$ explicitly, whence $(\xi(a)|\chi(a))$ too. We firstly obtain $(\xi(a)|\chi(a)) = A(x_3^3+x_4^2) + B$ where
\[\hspace{3 mm}  A ~ = ~ (q_3^2m_2 + q_1^2)q_1r_2^2b_3^{-2},\]
\[\begin{array}{lcl}
B  & = & (-m_1r_3 + (x_1-m_1b_2)q_1^2r_2b_3^{-1})x_1q_1^2b_3^{-1} + r_3^2m_2 +  \\
 & &  + (m_2-1)r_3r_2q_1^2b_2b_3^{-1} +  (x_1m_2q_3 - m_1q_3b_2)x_1q_3q_1^2r_2^2b_3^{-2}.
\end{array}\]
By exploiting (\ref{caso1 4}), (\ref{caso1 5}) and (\ref{caso12 3}) we eventually obtain the following:
\begin{equation}\label{caso12 6}
(\xi(a)|\chi(a)) = - r_3^2\frac{q_1^4m_2^2}{1+m_2} + q_3^2m_2 + q_1^2.
\end{equation} 
In this equation $(\xi(a)|\chi(a))$ is expressed as a function of $r_3$ rather than $b_3$, but recall that $r$ is uniquely determined by $b$. Note that the coefficient of $r_3^2$ in (\ref{caso12 6}) is negative except when $m_2 = 0$. If $m_2 = 0$ then $(\xi(a)|\chi(a)) = q_1^2$, which is strictly positive and less than $1$, since neither $q_1$ nor $q_3$ are zero. \\

\noindent
{\bf Case 2.} $\Gamma = \Gamma_\mathbb{C}(\mathbb{O})$. As in Case 1, we can assume that 
\[\begin{array}{lll}
L = [{\bf k},{\bf k}], & & M = [{\bf k}, {\bf j}m_1+{\bf k}m_2], ~~ |m_1|^2+|m_2|^2= 1,\\
p = [{\bf j}], & & q = [{\bf j}q_1 + {\bf kj}q_3], ~~ |q_1|^2 + |q_3|^2 = 1.
\end{array}\] 
So, $\ell(\alpha) = m_2$. As in Case 1, we have $q_1\neq 0 \neq q_3$. Let $N = [b,x]$ be admissible for $\alpha$, where 
\[\begin{array}{l}
b = {\bf j}b_1 + {\bf k}b_2 + {\bf kj}b_3, ~~ |b_1|^2 + |b_2|^2 + |b_3|^2 = 1,\\ 
x = {\bf j}x_1 + {\bf k}x_2 + {\bf kj}x_3, ~~ |x_1|^2 + |x_2|^2 + |x_3|^2 = 1.
\end{array}\]
For $N$ to be admissible for $\alpha$ the following must hold: $({\bf j}q_1+{\bf kj}q_3|b) = 0$ and $({\bf k}|b) = ({\bf k}|x) = ({\bf j}m_1+{\bf k}m_2|x)$. Explicitly: 
\begin{equation}\label{caso2 1}
\overline{q_1}b_1 + \overline{q_3}b_3 = 0, 
\end{equation}
and $b_2 = x_2 = \overline{m_1}x_1+\overline{m_2}x_2$, namely
\begin{equation}\label{caso2 2}
b_2 = x_2, \hspace{5 mm} \overline{m_1}x_1 = (1-\overline{m_2})b_2.
\end{equation}
Let $r = [{\bf j}r_1 + {\bf k}r_2 + {\bf kj}r_3]$ be such that $\{r\} = \{[b],p\}^\perp$. So, $r = [{\bf k}r_2 + {\bf kj}r_3]$, where we assume 
$|r_2|^2 + |r_3|^2 = 1$, and 
\begin{equation}\label{caso2 3}
\overline{r_2}b_2 + \overline{r_3}b_3 = 0.
\end{equation}
Recall that $q_1 \neq 0$ because $p\not\perp q$ by assumption. We also have $r_2 \neq 0$, otherwise $N$ cannot be coplanar with either of $L$ and $M$. Thus, by (\ref{caso2 1}) and 
(\ref{caso2 3}) we obtain 
\begin{equation}\label{caso2 4}
b_1 = -b_3\frac{\overline{q_3}}{\overline{q_1}}, \hspace{5 mm} b_2 = - b_3\frac{\overline{r_3}}{\overline{r_2}}.
\end{equation}
These equations show that $b_3 \neq 0$. Recalling that $|b_1|^2 + |b_2|^2 + |b_3|^2 = 1$ we get 
\begin{equation}\label{caso2 5}
b_3 ~ = ~ \varepsilon\cdot\frac{q_1r_2}{\sqrt{|q_1|^2 + |r_2|^2 - |q_1|^2|r_2|^2}} ~ = ~ \varepsilon\frac{q_1r_2}{\sqrt{1-|q_3|^2|r_3|^2}}
\end{equation}
for a suitable multiplier $\varepsilon$ with $|\varepsilon| = 1$. We shall now consider two subcases: either $|m_2| = 1$ or $|m_1| < 1$. \\

\noindent
{\bf Subcase 2.1.} $|m_2| = 1$. Equivalently, $m_1 = 0$. Then $b_2 = x_2 = 0$ by (\ref{caso2 2}), $r_3 = 0$ by (\ref{caso2 4}) and since $b_3 \neq 0$, whence $|r_2| = 1$ and $|b_3| = |q_1|$ by (\ref{caso2 5}). Consequently, $|b_1| = |q_3|$. 

Let now $\xi$ be the plane on $L$ and $N$ and $\chi$ the plane on $M$ and $N$. Then $\xi$ and $\chi$, regarded as  sharp $\mathbb{C}$-autorphisms of $\mathbb{O}$, are uniquely determined by the following conditions: $\xi({\bf k}) = {\bf k}$, $\chi({\bf k}) = {\bf j}m_1 + {\bf k}m_2$ and $\xi(b) = \chi(b) = x$. In view of the above:  
\begin{equation}\label{caso21 1}
\xi({\bf k}) = {\bf k}, ~~~ \chi({\bf k}) = {\bf k}m_2, ~~~ \xi(b) = \chi(b) = {\bf j}x_1 + {\bf kj}x_3.
\end{equation} 
It is easy to check that 
\[{\bf j} = ({\bf j}b_1 + {\bf kj}b_3)(\overline{b_1}+ {\bf k}\overline{b_3}) = b(\overline{b_1}+{\bf k}\overline{b_3}).\]
By this and (\ref{caso21 1}) we get
\begin{equation}\label{caso21 2}
\begin{array}{rcl}
\xi({\bf j}) & = & ({\bf j}x_1 + {\bf kj}x_3)(\overline{b_1} + {\bf k}\overline{b_3}),\\
\chi({\bf j}) & = & ({\bf j}x_1 + {\bf jkj}x_3)(\overline{b_1} + {\bf k}m_2\overline{b_3}).
\end{array}
\end{equation} 
Let $L' = [a,\xi(a)]$ and $M' = [a,\chi(a)]$ be the lines through $p$ and $r$ in $\xi$ and $\chi$ respectively, where 
$a = {\bf j}a_1 + {\bf k}a_2 + {\bf kj}a_3$  is orthogonal with both $p$ and $r$ and $|a_1|^2+|a_2|+|a_3|^2 = 1$. Orthogonality with $p$ and $r$ forces $a_1 = 0 = a_2$. Therefore $a = {\bf kj}\eta$ for a suitable $\eta$ with $|\eta| = 1$. By this and (\ref{caso21 2}), 
\begin{equation}\label{caso21 3}
\begin{array}{rcl}
\xi(a) & = & {\bf k}(({\bf j}x_1 + {\bf kj}x_3 + {\bf k})(\overline{b_1} + {\bf k}\overline{b_3}))\eta,\\
\chi(a) & = & {\bf k}m_2(({\bf j}x_1 + {\bf kj}x_3 + {\bf k}x_4)(\overline{b_1} + {\bf k}m_2\overline{b_3}))\eta.
\end{array}
\end{equation} 
Equations (\ref{caso21 3}) allow to explicitly compute the inner product $(\xi(a)|\chi(a))$. We obtain:
\begin{equation}\label{caso21 4}
(\xi(a)|\chi(a)) ~=~ |b_3|^2 + |b_1|^2\overline{m_2} ~ = ~ |q_1|^2 + |q_3|^2\overline{m_2}.
\end{equation}
So, $|(\xi(a)|\chi(a))| = |q_1|^4+|q_3|^4 + |q_1|^2|q_3|^2(m_2+\overline{m_2}) < 1$, as $m_2+\overline{m_2}$ is a real number not less than $-2$ and less than $2$ (because $|m_2| = 1$ but $m_2\neq 1$) and $|q_1|^2+|q_3|^2 = 1$. \\

\noindent
{\bf Subcase 2.2.} $m_1 \neq 0$, namely $|m_2| < 1$. In this case the second equation of (\ref{caso2 2}) yields
\begin{equation}\label{caso22 1}
x_1 ~ = ~ \frac{1-\overline{m_2}}{\overline{m_1}}b_2.
\end{equation} 
The planes $\xi$ and $\chi$ on $L$ and $N$ and on $M$ and $N$ are determined by the following conditions:
\begin{equation}\label{caso22 2}
\begin{array}{l}
\xi({\bf k}) = {\bf k}, \hspace{5 mm} \chi({\bf k}) = {\bf j}m_1 + {\bf k}m_2, \\
\xi(b) = \chi(b)  =  {\bf j}x_1 + {\bf k}x_2 + {\bf kj}x_3  =  ({\bf j}\frac{1-m_2}{m_1} + {\bf k})b_2 + {\bf kj}x_3.
\end{array}
\end{equation} 
Moreover, $|x_3|^2 = 1- (1+|1-m_2|^2|m_1|^{-2})|b_2|^2$ by (\ref{caso22 1}) and $x_2 = b_2$. Therefore
\begin{equation}\label{caso22 3}
|x_3|^2 ~ = ~ 1 - \frac{2-m_2-\overline{m_2}}{|m_1|^2}|b_2|^2.
\end{equation}
Now ${\bf j} = (b-{\bf k}b_2)((\overline{b_1}+{\bf k}\overline{b_3})(1-|b_2|^2)^{-1})$. Recalling equations (\ref{caso2 4}), we obtain
\begin{equation}\label{caso22 4}
{\bf j} ~ = ~ (b+{\bf k}\frac{\overline{r_3}}{\overline{r_2}}b_3)(({\bf k}q_1 - q_3)\overline{q_1}b_3^{-1}).
\end{equation} 
Let $L' = [a,\xi(a)]$ and $M' = [a,\chi(a)]$ be the lines through $p$ and $r$ in $\xi$ and $\chi$ respectively, where 
$a = {\bf j}a_1 + {\bf k}a_2 + {\bf kj}a_3$ is orthogonal with both $p$ and $r$ and $|a| = 1$. Orthogonality with $p$ forces $a_1 = 0$ but orthogonality with $r$ only implies $\overline{r}_2a_2 + \overline{r_3}a_3 = 0$. So $a_2 = - a_3\overline{r_3}\overline{r_2}^{-1}$ and the condition $|a| = 1$ implies $|a_3| =  |r_2|$, namely $a_3 = \overline{r_2}\eta$ for some $\eta$ with $|\eta| = 1$. Hence
\begin{equation}\label{caso22 5}
a = (-{\bf k}\overline{r_3} + {\bf kj}\overline{r_2})\eta = ({\bf k}(-\overline{r_3} + \overline{r_2}{\bf j}))\eta =  ({\bf k}(-\overline{r_3} + {\bf j}r_2))\eta.
\end{equation}
By exploiting (\ref{caso22 2}), (\ref{caso22 4}) and (\ref{caso22 5}) as well as (\ref{caso2 4}) and (\ref{caso22 3}) one can compute $\xi(a)$ and $\chi(a)$ explicitly, whence $(\xi(a)|\chi(a))$ too, but these computations are terribly toilsome. However, in order to prove the lemma, we do not need to perform them.  
It is enough to show that, for a lucky choice of $N = [b,x]$, whence of $r$, satisfying the above conditions, we get $\ell \neq -1$. We will go on in this way, referring the interested reader to Remark \ref{polinomiaccio} for a way to express $(\xi(a)|\chi(a))$ in the general case.

The previous conditions on $r$, $b$ and $x$ allow to choose $r_3 = 0$. Accordingly, $|r_2| = 1$. Hence $b_2 = 0$ by the second equation of (\ref{caso2 4}) and $b_3 = \lambda\overline{q_1}$ for some $\lambda$ with $|\lambda| = 1$ by (\ref{caso2 5}). Therefore $b_1 = -\lambda\overline{q}_3$ by the first equation of (\ref{caso2 4}).  Moreover $x_1 = x_2 = 0$ by (\ref{caso2 2}) and (\ref{caso22 1}), whence $|x_3| = 1$. Accordingly, 
\begin{equation}\label{caso22 4 bis}
{\bf j} = b(({\bf k}q_1-q_3)\lambda^{-1})
\end{equation}
by (\ref{caso22 4}) and since $b_1 = \lambda\overline{q}_1$ and
\begin{equation}\label{caso22 5 bis}
a = {\bf kj}\overline{r_2}\eta
\end{equation}
by (\ref{caso22 5}) and since $r_3 = 0$. By (\ref{caso22 4 bis}), recalling that $x_1 = x_2 = 0$, we obtain 
\begin{equation}\label{caso22 6}
\left.\begin{array}{rcl}
\xi({\bf j}) & = &  x(({\bf k}q_1-q_3)\lambda^{-1}) = {\bf kj}x_3({\bf k}q_1\lambda^{-1} -q_3\lambda^{-1}) =\\
 & = & {\bf j}\overline{q_1}\overline{x_3}\lambda - {\bf kj}q_3x_3\overline{\lambda},\\
\chi({\bf j}) & = &  x((({\bf j}m_1+{\bf k}m_2)q_1-q_3)\lambda^{-1}) = \\
 & = & {\bf kj}x_3({\bf j}m_1q_1\lambda^{-1} + {\bf k}m_2q_1\lambda^{-1}-q_3\lambda^{-1}) = \\
 & = & {\bf j}\overline{m_2}\overline{q_1}\overline{x_3}\lambda -{\bf k}\overline{m_1}\overline{q_1}\overline{x_3}\lambda - {\bf kj}q_3x_3\overline{\lambda}.
\end{array}\right\}
\end{equation} 
(Recall that $\lambda^{-1} = \overline{\lambda}$ since $|\lambda| = 1$.) By combining (\ref{caso22 5 bis}) with (\ref{caso22 6}) we obtain
\[\begin{array}{rcl}
\xi(a) & = & ({\bf k}({\bf j}\overline{q_1}\overline{x_3}\lambda - {\bf kj}q_3x_3\overline{\lambda}))\overline{r_2}\eta = \\
 & = & {\bf j}\overline{q_3}\overline{x_3}\overline{r_2}\lambda\eta + {\bf kj}q_1x_3\overline{r_2}\overline{\lambda}\eta,\\
\chi(a) & = & (({\bf j}m_1+{\bf k}m_2)({\bf j}\overline{m_2}\overline{q_1}\overline{x_3}\lambda -{\bf k}\overline{m_1}\overline{q_1}\overline{x_3}\lambda - {\bf kj}q_3x_3\overline{\lambda}))\overline{r_2}\eta = \\
 & = & {\bf j}\overline{m_2}\overline{q_3}\overline{x_3}\overline{r_2}\lambda\eta - {\bf k}\overline{m_1}\overline{q_3}\overline{x_3}\overline{r_2}\lambda\eta + {\bf kj}q_1x_3\overline{r_2}\overline{\lambda}\eta.
\end{array}\]
Therefore $(\xi(a)|\chi(a)) =  (|q_3|^2\overline{m_2} + |q_1^2)(|x_3|^2|r_2|^2|\lambda|^2|\eta|^2$. Finally, recalling that $|x_3| = |r_2| = |\lambda| = |\eta| = 1$,
\begin{equation}\label{caso22 7} 
(\xi(a)|\chi(a)) ~ = ~  |q_3|^2\overline{m_2}+|q_1|^2.
\end{equation}
The right side of (\ref{caso22 7}) is equal to $-1$ only if $q_1 = 0$ and $m_2 = -1$. However, $q_1 \neq 0$ because $p\not\perp q$. Therefore $(\xi(a)|\chi(a)) \neq -1$.  \hfill $\Box$

\begin{note}\label{polinomiaccio}
{\rm In Subcase 2.2 of the above proof, with no additional hypotheses on $[b,x]$ we get}
\[(\xi(a)|\chi(a)) ~=~ |q_3|^2\overline{m_2} + |q_1|^2 - 2\mathrm{Im}(m_1\overline{q_1}q_3|q_3|^2r_2\overline{r_3}x_3b_3^{-1}) + |r_3|^2A\]
{\rm where $\mathrm{Im}(.)$ stands for imaginary part and} 
\[\begin{array}{rcl}
A & :=  & B + m_2 - (|q_3|^2\overline{m_2}+|q_1|^2)(|q_3|^2 + |q_1|^2(2-m_2-\overline{m_2})|m_1|^{-2}),\\
B & := & |q_1|^2\overline{m_2}((1-m_2)^2|q_1|^2 - m_2(1-\overline{m_2})^2|q_3|^2)|m_1|^{-2}.
\end{array} \]
{\rm This shows that $(\xi(a)|\chi(b))$ depends on $r_2, r_3$ and $x_2$ non-trivially. Thus, we can always choose the line $N = [b,x]$ in such a way that $|(\xi(a)|\chi(a))| < 1$. Accordingly, Lemma \ref{ortho-invariant} can be given a stronger formulation: we can always choose $N$ in such a way that $|\ell| < 1$.}
\end{note}

\begin{note}
{\rm It follows from above proof that when $|m_2| = 1$ then $|\ell| < 1$ for every choice of the admissible line $N = [b,x]$. However, for certain values of $m_2$ we can also choose $N$ in such a way that $\ell = -1$. For instance, when $(\mathbb{F},\mathbb{A}) = (\mathbb{R},\mathbb{H})$, this is possible in the following cases:

\begin{itemize}
\item[1)] $q_1^4 = q_3^2$ (namely $q_1^2 = (\sqrt{5}-1)/2$) and $-1\leq m_2 \leq -(\sqrt{5}+1)/4$;
\item[2)] $q_1^2 > q_3^2$ and $-1\leq m_2 \leq (1-\sqrt{4q_1^6+8q_1^4-3})/(q_1^4-q_3^2)$;
\item[3)] $q_1^2 < q_3^2$ and $1\geq m_2 \geq (-1+\sqrt{4q_1^6+8q_1^2-3})/(q_3^2-q_1^4)$.
\end{itemize}} 
\end{note}

\begin{lemma}\label{SS 6.7}
Every orthogonal non-degenerate primitive path $\alpha$ of $\Gamma_\mathbb{C}(\mathbb{O})$ such that $|\ell(\alpha)| = 1$ but $\ell(\alpha) \neq -1$ is homotopic with an orthogonal non-degenerate primitive path $\beta$ such that $|\ell(\beta)| < 1$.
\end{lemma}

See Schillewaert and Struyve \cite[Lemma 6.7]{SS} for the above. The following lemma is also proved in \cite[Lemma 6.8]{SS}.

\begin{lemma}\label{SS 6.8}
Let $\ell\in \mathbb{F}$ such that $|\ell| < 1$. Then, fort every choice of two distinct lines $L$ and $M$ with the same shadow, there exists a sequence $L_0 = L, L_1,..., L_n = M$ of lines with the same sadow as $L$ and $M$ and such that $(L_{i-1}|L_i) = \ell$ for every $i = 1, 2,..., n$. 
\end{lemma}

The next statement is implicit in what Schillewaert and Struyve say to justify \cite[Remark 6.9]{SS}. We make it explicit. 

\begin{cor}\label{SS 6.8 cor}
Let $\ell\in \mathbb{F}$ such that $|\ell| < 1$ and let $\alpha = (p,L,q,M,p)$ be a non-degenerate primitive path of $\Gamma = \Gamma_\mathbb{F}(\mathbb{A})$. Then $\alpha \sim \alpha_1\cdot\alpha_2\cdot...\cdot\alpha_n$ for a suitable sequence of non-degnenerate primitive paths 
$\alpha_1, \alpha_2,..., \alpha_n$ of $\Gamma$ with the same points $p$ and $q$ as $\alpha$ and such that $\ell(\alpha_i) = \ell$ for every $i = 1, 2,..., n$.
\end{cor}
{\bf Proof.} By Lemma \ref{SS 6.8} there exist lines $L_0 = L, L_1,..., L_n = M$ such that $(L_{i-1}|L_i) = \ell$ for $i = 1, 2,..., n$. For $i = 1, 2,..., n$ put $\alpha_i = (p,L_{i-1},q,L_i)$. Thus, the concatenation $\alpha_1\cdot\alpha_2\cdot...\cdot\alpha_n$ is well defined. Note that
\[\alpha_{n-1}\cdot\alpha_n = (p, L_{n-2}, q, L_{n-1}, p, L_{n-1}, q, L_n, p) \sim (p, L_{n-2}, q, L_n) = : \alpha'_{n-1}.\]
So, $\alpha_1\cdot\alpha_2\cdot...\cdot\alpha_{n-1}\cdot\alpha_n \sim \alpha_1\cdot\alpha_3\cdot...\cdot\alpha'_{n-1}$. By iterating this argument we eventually obtain $\alpha_1\cdot\alpha_2\cdot...\cdot\alpha_n \sim (p, L_0, q, L_n, p) = \alpha$. \hfill $\Box$ \\

We can now prove the main theorem of this subsection.

\begin{theo}\label{theo2}
Either $\Gamma_\mathbb{F}(\mathbb{A})$ is simply connected or it is covered by a building.
\end{theo}
{\bf Proof.} Suppose that $\Gamma = \Gamma_\mathbb{F}(\mathbb{A})$ is not covered by a building. Then, by Corollary \ref{primitive cor2}, at least one of its non-degenerate primitive paths is null-homotopic. By Lemma \ref{shift} (claim (4)) and Lemma \ref{ortho-shift}, at least one orthogonal non-degenerate primitive path, say $\alpha$, is null-homotopic. Let $\ell = \ell(\alpha)$ be its line-invariant. The action of $G := \mathrm{Aut}(\Gamma)$ on $\mathbb{A}$ and $\mathbb{O}$ makes it clear that $G$ acts transitively on the set of orthogonal primitive paths with line-invariant equal to $\ell$. Thus, all orthogonal primitive paths with line invariant $\ell$ are null-homotopic. 

Suppose firstly that $|\ell| < 1$. Then every orthogonal primitive path $\beta$ is null homotopic, by Corollary \ref{SS 6.8 cor} and the above remark. In this case $\Gamma$ is simply connected by Lemmas \ref{shift} and \ref{ortho-shift} and Corollary \ref{reduction3}.

Let $|\ell| = 1$. If $\ell \neq -1$ (whence $\Gamma = \Gamma_\mathbb{C}(\mathbb{O})$) then $\alpha \sim \beta$ for some orthogonal primitive path $\beta$ with $|\ell(\beta)| < 1$, by Lemma \ref{SS 6.7}. Thus, we can replace $\alpha$ with $\beta$ and we are back to the previous case. 

Finally, let $\ell(\alpha) = -1$. Clearly $\alpha$ admits a non-orthogonal shift $\beta \sim \alpha$, necessarily non-degenerate (Lemma \ref{shift}). In its turn $\beta$ admits an orthogonal shift $\gamma$ with $\ell(\gamma) \neq -1$, by Lemma \ref{ortho-invariant}. Moreover $\beta\sim \gamma$ by Lemma \ref{shift}. Hence $\alpha\sim \gamma$. Therefore $\gamma$ is null-homotopic. We can now replace $\alpha$ with $\gamma$ and we are back to the first or second one of the two previous cases, according to whether $|\ell(\gamma)| < 1$ or $|\ell(\gamma)| = 1$.  \hfill $\Box$  

\begin{note}
{\rm What Schillewaert and Struyve say to explain \cite[Remark 6.9]{SS} amounts to a sketch of the first three paragraphs of the above proof but, as they had nothing like Lemma \ref{ortho-invariant} at their disposal, they necessarily missed the very last step of the proof.}
\end{note}       

\subsection{End of the proof of Theorem \ref{prel3}}

Let $\Gamma = \Gamma_\mathbb{F}(\mathbb{A})$ and let $\phi_\Gamma:\tGamma\rightarrow\Gamma$ be its universal covering. In view of Theorem \ref{theo2}, either $\tGamma = \Gamma$ or $\tGamma$ is a building. In order to finish the proof of Theorem \ref{prel3} it only remains to prove that $\tGamma$ cannot be a building. This immediately follows from Theorem \ref{prel2}. However, as we have promised not to use that theorem, we shall give an explicit proof of this claim. We firstly state some notation.   

\subsubsection{Some notation}

For a positive integer $n$, let $f^\mathbb{F}_n$ be the usual scalar product on $\mathbb{F}^n$ and let $\mathrm{L}(f^\mathbb{F}_n)$ be the group of all linear mappings preserving $f^\mathbb{F}_n$. So, $\mathrm{L}(f^\mathbb{R}_n) = \mathrm{O}(n)$ and $\mathrm{L}(f^\mathbb{C}_n) = \mathrm{U}(n)$ (notation as usual for Lie groups). 

Given two positive integers $n, m$ with $n\leq m$, let $f^\mathbb{F}_{n,m} := (-f^\mathbb{F}_n)\oplus f^\mathbb{F}_m$. Namely, $f^\mathbb{F}_{n,m}$ admits the following representations, according to whether $\mathbb{F} = \mathbb{R}$ or $\mathbb{F} = \mathbb{C}$, where $x = (x_i)_{i=1}^{n+m}$ and $y = (y_i)_{i=1}^{n+m}$ (vectors of $\mathbb{F}^{n+m}$):
\[\begin{array}{lcrcl}
(\mathbb{F} = \mathbb{R}) & & f_{n,m}^\mathbb{R}(x,y) & := &  -\sum_{i=1}^nx_iy_i  + \sum_{i=1}^mx_{i+n}y_{i+m},\\
(\mathbb{F} = \mathbb{C}) & & f_{n.m}^\mathbb{C}(x,y) & := &  -\sum_{i=1}^n\overline{x_i}y_i  + \sum_{i=1}^m\overline{x_{i+n}}y_{i+m}.
\end{array}\]
Clearly, $n$ is the Witt index of $f^\mathbb{F}_{n,m}$. We also recall that every non-degenerate bilinear form on $\mathbb{R}^{n+m}$ of Witt index $n \leq m$ can be expressed as $f^\mathbb{R}_{n,m}$ or its opposite, modulo a suitable choice of the basis of $\mathbb{R}^{n+m}$. The same for hermitian forms of $\mathbb{C}^{n+m}$. 

Let $\mathrm{L}(f^\mathbb{F}_{n.m})$ be the group of linear trasformations of $\mathbb{F}^{n+m}$ preserving $f^\mathbb{F}_{n,m}$. So, $\mathrm{L}(f^\mathbb{R}_{n,n}) = \mathrm{O}^+(2n,\mathbb{R})$, $\mathrm{L}(f^\mathbb{R}_{n,n+1}) = \mathrm{O}(2n+1, \mathbb{R})$, $\mathrm{L}(f^\mathbb{C}_{n,n}) = U(2n, \mathbb{C})$ and $\mathrm{L}(f^\mathbb{C}_{n,n+1}) = U(2n+1, \mathbb{C})$ (notation as usual for Chevalley groups). 

Let $\Gamma(f^\mathbb{F}_{n,m})$ be the polar space associated to $f^\mathbb{F}_{n,m}$. Recall that its full automorphims group $\mathrm{Aut}(\Gamma(f^\mathbb{F}_{n,m}))$ is the projectivization $\mathrm{PL}(f^\mathbb{F}_{n,m})$ of $\mathrm{L}(f^\mathbb{F}_{n,m})$, extended by two (possibly trivial) outer automorphism groups, henceforth denoted by $\bf d$ and $\bf f$. The group $\bf d$ is contributed by linear transformations of $\mathbb{F}^{n+m}$ which do not preserve $f^\mathbb{F}_{n,m}$ but multiply it by a scalar. However, as we deal with $\mathrm{PL}(f^\mathbb{F}_{n,m})$ rather than $\mathrm{L}(f^\mathbb{F}_{n,m})$, it turns ut that $\bf d$ is either trivial or isomorphic to the group $C_2$ of order $2$, according to whether $n+m$ is odd or even. The group $\bf f$ is trivial when $\mathbb{F} = \mathbb{R}$ and isomorphic to $C_2$ when $\mathbb{F} = \mathbb{C}$. In the latter case, the unique non trivial involution of $\bf f$ is contributed by the usual conjugation of $\mathbb{C}$ and the extension $(\mathrm{PL}(f^\mathbb{C}_{n,m})\cdot{\bf d})\cdot{\bf f}$ is split: it can be realized as the semidirect product $(\mathrm{PL}(f^\mathbb{C}_{n,m})\cdot{\bf d})>\hspace{-2.4mm}\lhd \langle\iota\rangle$ of $\mathrm{PL}(f^\mathbb{C}_{n,m})\cdot {\bf d}$ with the group $\langle \iota\rangle$ generated by a suitable involutory semi-linear transformation $\iota$ of $\mathbb{C}^{n+m}$.  

\subsubsection{The case $(\mathbb{F}, \mathbb{A}) = (\mathbb{C}, \mathbb{O})$}\label{fine 1}

Let $\Gamma = \Gamma_\mathbb{C}(\mathbb{O})$. By contradiction, suppose that $\tGamma$ is a building. Then, by considering dimensions of panels, we see that $\tGamma = \Gamma(f^\mathbb{C}_{3,4})$, with full automorphism group 
$\widehat{G} := \mathrm{Aut}(\Gamma(f^\mathbb{C}_{3,4}))  =  \mathrm{PU}(7,\mathbb{C})>\hspace{-2.4mm}\lhd  {\bf f}  =  \mathrm{PSU}(7,\mathbb{C})>\hspace{-2.4mm}\lhd   C_2$.

We set $G := \mathrm{Aut}(\Gamma) = ((\mathrm{SU}(3)\times\mathrm{SU}(3))/C_3)>\hspace{-2.4mm}\lhd   C_2$ (see Section \ref{Aut groups}). 

Let $\tilde{\xi}$ be a plane of $\tGamma$ and $\xi = \phi_\Gamma(\tilde{\xi})$. Let $G_\xi$ be the stabilizer of $\xi$ in $G$ and $\widehat{G}_{\tilde{\xi}}$ the stabilizer of $\tilde{\xi}$ in $\widehat{G}$. The group $G_\xi$ should be recognizable as a subgroup of $\widehat{G}_{\tilde{\xi}}$. More explicitly: $G_\xi \cong \tG_\xi$ for a suitable subgroup $\tG_\xi < \widehat{G}_{\tilde{\xi}}$. It is not so difficult to see that $G_\xi  =  \mathrm{PSU}(3)>\hspace{-2.4mm}\lhd  C_2$ $(\cong \tG_\xi)$. On the other hand, $\widehat{G}_{\tilde{\xi}} = U>\hspace{-2.4mm}\lhd  L$ where 
$L \cong \mathrm{GL}(3,\mathbb{C})>\hspace{-2.4mm}\lhd {\bf f} = \Gamma\mathrm{L}(3,\mathbb{C})$ and $U \cong \mathbb{C}^6\times\mathbb{R}^3 ~ \cong ~ \mathbb{R}^{15}$, with $\mathbb{C}^6$, $\mathbb{R}^3$ and $\mathbb{R}^{15}$ being regarded as additive groups. Needless to say, $U$ is the unipotent radical of $\widehat{G}_{\tilde{\xi}}$ and $L$ plays the role of Levi complement. 

We have $\tG_\xi\cap U = 1$, since $\mathrm{PSU}(3)>\hspace{-2.4mm}\lhd  C_2$ admits no infinite commutative normal subgroups. Hence $\tG_\xi \leq L \cong \Gamma\mathrm{L}(3,\mathbb{C})$. The group $\Gamma\mathrm{L}(3,\mathbb{C})$ indeed contains copies of $\mathrm{SU}(3)>\hspace{-2.4mm}\lhd  C_2$, but no copy of  $\mathrm{PSU}(3)>\hspace{-2.4mm}\lhd  C_2$. We have reached a contradiction. Hence in this case $\tGamma = \Gamma$. 

\subsubsection{The case $(\mathbb{F},\mathbb{A}) = (\mathbb{R},\mathbb{H})$}\label{fine 2}

Let $\Gamma = \Gamma_\mathbb{R}(\mathbb{H})$. By contradiction, suppose that $\tGamma$ is a building. Then, by considering dimensions of panels, we see that $\tGamma = \Gamma(f^\mathbb{R}_{3,8})$, with full automorphism group $\widehat{G} := \mathrm{Aut}(\Gamma(f^\mathbb{R}_{3,8})) = \mathrm{L}(f^\mathbb{R}_{3,8})$. We set $G := \mathrm{Aut}(\Gamma) = \mathrm{SO}(3)\times \mathrm{G}_2$ (see Section \ref{Aut groups}). 

As above, let $\tilde{\xi}$ be a plane of $\tGamma$, let $\xi = \phi_\Gamma(\tilde{\xi})$ and let $G_\xi$ be the stabilizer of $\xi$ in $G$ and $\widehat{G}_{\tilde{\xi}}$ the stabilizer of $\tilde{\xi}$ in $\widehat{G}$. The group $G_\xi$ should be recognizable as a subgroup of $\widehat{G}_{\tilde{\xi}}$. It is not so difficult to check that
\[ G_\xi ~ = ~ (\mathrm{SU}(2)\times \mathrm{SU}(2))/\langle (-\iota,-\iota)\rangle ~ = ~ 2^\cdot (\mathrm{PSU}(2)\times\mathrm{PSU}(2)).\]
Here $\iota$ stands for the identity element of $\mathrm{SU}(2)$, whence $(\iota,\iota)$ is the identity element of $\mathrm{SU}(2)\times \mathrm{SU}(2)$. The extension $2^\cdot (\mathrm{PSU}(2)\times\mathrm{PSU}(2))$ is non-split. On the other hand, $\widehat{G}_{\tilde{\xi}} =  U>\hspace{-2.4mm}\lhd  L$ where 
\[L ~ \cong ~ (\mathrm{GL}(3,\mathbb{R})\times \mathrm{O}(5))/\langle(-\iota,-\iota)\rangle ~ \cong ~ \mathrm{GL}(3,\mathbb{R})\times \mathrm{SO}(5)\]
and $U = {U_0}^\cdot U_1$ with $U_0 \cong \mathrm{R}^3$ and $U_1\cong \mathrm{R}^{15}$. The extension ${U_0}^\cdot U_1$ is non-split. Explicitly, if $u\in U\setminus U_0$ then $1 \neq u^2\in U_0$. Therefore, every non-trivial subgroup of $U$ contains a non trivial subgroup of $U_0$, which is clearly commutative and infinite. The group $G_\xi$ contains no normal commutative infinite subgroup. Hence $\tG_\xi\cap U = 1$. Consequently, $\tG_\xi  \leq  L  \cong  \mathrm{GL}(3,\mathbb{R})\times \mathrm{SO}(5)$. However, $\tG_\xi \cong G_\xi$ is a non-split extension $2^\cdot (\mathrm{PSU}(2)\times\mathrm{PSU}(2))$. No group with this structure can be hosted as a subgroup by $\mathrm{GL}(3,\mathbb{R})\times \mathrm{SO}(5)$. Again, a contradiction. 

Therefore $\tGamma = \Gamma$ in this case too. The proof of Theorem \ref{prel3} is complete. 

\begin{note}
{\rm In order to prove that $\Gamma$ is not covered by a building, Kramer and Lytchak \cite{KL1} and \cite{KL2} also describe the structure of a chamber-stabilizer in $G$, showing that, in either of the two possible cases, it does not fit with any subgroup of a chamber-stabilizer in $\widehat{G}$. The arguments they use to prove this fact are not so different from those we have exploited here but more toilsome, due to the fact that the structure of a chamber-stabilizer in either $G$ or $\widehat{G}$ is fairly poor if compared with that of a plane-stabilizer.}
\end{note}

\bigskip

\noindent
Author's address\\

\noindent
Antonio Pasini,\\ 
retired,\\
antonio.pasini@unisi.it


\begin{thebibliography}{99}
\bibitem{Asch} M. Aschbacher, {\em Finite geometries of type $C_3$ with flag-transitive automorphism groups}, Geom. Dedicata {\bf 16} (1984), 195-200.
\bibitem{BC} A. E. Brouwer and A. M. Cohen, {\em Some remarks on Tits's geometries}, Indagationes Math. {\bf 45} (1983), 393-402.  
\bibitem{BP} F. Buekenhout and A. Pasini, {\em Finite diagram geometries extending buildings}, Chapter 22 of Handook of Incidence Geometry (ed. F. Buekenhout), North-Holland, Amsterdam (1995), 1143-1254.  
\bibitem{GK} G. Goroski and A. Kollross, {\em Some remarks on polar actions}, Ann. Global Anal. Geom. {\bf 49} (2016), 43-58.  
\bibitem{KL1} L. Kramer and A. Lytchak, {\em Homogeneous compact geometries}, Transform. Groups {\bf 19} (2014), 793-852.
\bibitem{KL2} L. Kramer and A. Lytchak, {\em Erratum to: Homogeneous compact geometries}, Transform. Groups, to appear. 
\bibitem{N} A. Neumaier, {\em Some sporadic geometries related to $PG(3,2)$}, Arch. Math. {\bf 42} (1984), 89-96. 
\bibitem{PasDG} A. Pasini, Diagram Geometries, Oxford Univ. Press, 1994. 
\bibitem{PT} F. Podest\`{a} and G. Thorbergsson, {\em Polar actions on rank-one symmetric spaces}, J. Differential Geom. {\bf 53} (1999), 131-175. 
\bibitem{Rees} S. E. Rees, {\em $C_3$ geometries arising from the Klein quadric}, Geom. Dedicata {\bf 18} (1985), 67-85. 
\bibitem{SS} J. Schillewaert  and K. Struyve, {\em On exceptional homogeneous compact geometries of type $C_3$}, Groups Geom. Dyn. {\bf 11} (2017), 1377-1399.  
\bibitem{Tits1} J. Tits, Buildings of Spherical Type and Finite $BN$-Pairs, Springer L.N. {\bf 386}, Springer, Berlin, 1974. 
\bibitem{Tits2} J. Tits, {\em A local approach to buildings}, in The Geometric Vein (eds. C. Davies et al.), Springer, Berlin (1981), 519-547.  
\bibitem{Yoshi} S. Yoshiara, {\em Flag-transitive $C_3$-geometries of finite order}, J. Alg. Combin. {\bf 5} (1996), 251-284. 
\end{thebibliography}
\end{document}